\documentclass[aos]{imsart}
\usepackage{xr}
\RequirePackage{amsthm,amsmath,amsfonts,amssymb}

\RequirePackage{natbib}

\RequirePackage[colorlinks,citecolor=blue,urlcolor=blue]{hyperref}

\usepackage{verbatim}
\usepackage{natbib,graphicx, lscape,longtable}
\usepackage{natbib,epsfig,graphicx}
\usepackage{booktabs}
\usepackage{rotating}
\usepackage{multirow}
\usepackage{amsmath,amsthm,amssymb,color}
\bibpunct{(}{)}{;}{a}{,}{,}

\newtheorem{tm}{Theorem}
\newtheorem{as}{Assumption}

\newtheorem{re}{Remark}
\newtheorem{pn}{Proposition}

\def\mF{\mathcal F}

\def\mG{\mathcal G}

\def\mB{\mathcal B}
\def\mV{\mathcal V}
\def\tr{{\rm tr}}
\def\mX{\mathcal X}

\def\mE{\mathcal E}

\def\bmC{\mathbf{C}}

\def\be{\begin{equation}}
\def\ee{\end{equation}}
\def\ben{\begin{equation*}}
\def\een{\end{equation*}}
\def\bea{\begin{eqnarray}}
\def\eea{\end{eqnarray}}

\def\bda{\begin{eqnarray*}}
\def\eda{\end{eqnarray*}}

\def\lsk{\left(}
\def\rsk{\right)}
\def\lbk{\left \{ }
\def\rbk{\right \} }
\def\lmk{\left [ }
\def\rmk{\right ] }
\def\lak{\left | }
\def\rak{\right | }
\def\laak{\left \| }
\def\raak{\right \| }

\def\nn{\nonumber}

\def\bmC{\mathbf{C}}
\def\bmV{\mathbf{V}}

\def\bmX{\mathbf{X}}
\def\bmY{\mathbf{Y}}

\def\bmW{\mathbf{W}}
\def\mbW{\mathbb{W}}
\def\mbV{\mathbb{V}}

\def\bmx{\mathbf{x}}

\def\rmE{{\rm E}}
\def\rmP{{\rm P}}

\def\rmvec{{\rm vec}}
\def\Cov{{\rm Cov}}

\def\tr{{\rm tr}}
\def\Cov{{\rm Cov}}

\def\mbX{\mathbb{X}}
\def\mbY{\mathbb{Y}}

\def\ErrF{{\rm ErrF}}
\def\ErrS{{\rm ErrS}}
\def\ErrR{{\rm ErrR}}
\def\Err{{\rm Err}}

\def\OptF{{\rm OptF}}
\def\OptS{{\rm OptS}}
\def\OptR{{\rm OptR}}
\def\Opt{{\rm Opt}}
\def\Avar{{\rm Avar}}
\def\QMLE{{\rm QMLE}}
\def\OLS{{\rm OLS}}
\def\Cp{{\rm Cp}}
\def\RCp{{\rm RCp}}
\def\WLS{{{\rm WLS}}}
\def\rmG{{{\rm W}}}
\def\GLS{{{\rm GLS}}}
\def\FGLS{{\rm FGLS}}
\def\Tr{{\rm Tr}}
\def\Te{{\rm Te}}
\def\OCV{{\rm OCV}}
\def\wbmC{{\widetilde \bmC}}

\numberwithin{equation}{section}

\newcommand{\lrnormmm}[1]{{\left\vert\kern-0.25ex\left\vert\kern-0.25ex\left\vert #1 \right\vert\kern-0.25ex\right\vert\kern-0.25ex\right\vert}}
\newcommand{\normmm}[1]{{\vert\kern-0.25ex\vert\kern-0.25ex\vert #1 \vert\kern-0.25ex\vert\kern-0.25ex\vert}}

\newcommand{\lrnormm}[1]{{\left\vert\kern-0.25ex\left\vert #1 \right\vert\kern-0.25ex\right\vert}}
\newcommand{\normm}[1]{{\vert\kern-0.25ex\vert #1 \vert\kern-0.25ex\vert}}

\def\singlespace{\def\baselinestretch{1}\@normalsize}

\begin{document}

\begin{frontmatter}
\title{Fixed and Random Covariance Regression Analyses}
\runtitle{Fixed and Random Covariance Regression}

\begin{aug}
\author[A]{\fnms{Tao}~\snm{Zou}\ead[label=e1]{tao.zou@anu.edu.au}},
\author[B]{\fnms{Wei}~\snm{Lan}\ead[label=e2]{lanwei@swufe.edu.cn}},

\author[C]{\fnms{Runze}~\snm{Li}\ead[label=e3]{rzli@psu.edu}}
\and
\author[D]{\fnms{Chih-Ling}~\snm{Tsai}\ead[label=e4]{cltsai@ucdavis.edu}}
\address[A]{The Australian National University\printead[presep={,\ }]{e1}}
\address[B]{Southwestern University of
Finance and Economics\printead[presep={,\ }]{e2}}
\address[C]{The Pennsylvania State University\printead[presep={,\ }]{e3}}
\address[D]{University of California, Davis\printead[presep={,\ }]{e4}}
\end{aug}

\begin{abstract}
Covariance regression analysis is an approach to linking the covariance of responses to a  set of explanatory variables $X$, where $X$ can be a vector, matrix, or  tensor. Most of the literature on this topic focuses on the ``Fixed-$X$'' setting and treats $X$ as nonrandom. By contrast, treating explanatory variables $X$ as  random, namely the ``Random-$X$'' setting, is often more realistic in practice. This article aims to fill this gap in the literature on the estimation and  model assessment theory for Random-$X$ covariance regression models. Specifically, we construct a new theoretical framework for studying the covariance estimators under the Random-$X$ setting, and we demonstrate that the quasi-maximum likelihood estimator and the weighted least squares estimator are both consistent and asymptotically normal. In addition, we develop pioneering work on the  model assessment theory of covariance regression. In particular,  we obtain the bias-variance decompositions for the expected test errors under both the Fixed-$X$ and Random-$X$  settings.
We show that moving from a Fixed-$X$ to a Random-$X$ setting can increase both the bias and the variance in  expected test errors.
 Subsequently, we propose estimators of the  expected  test errors under the Fixed-$X$ and Random-$X$  settings, which can be used to assess the performance of the competing covariance regression models. The proposed estimation and model assessment approaches are illustrated via extensive simulation experiments and an empirical study of stock returns in the US market.
\end{abstract}

\begin{keyword}[class=MSC]
\kwd[Primary ]{62H12}
\kwd[; secondary ]{62J02}
\end{keyword}

\begin{keyword}
\kwd{Bias-variance decomposition}
\kwd{covariance regression}
\kwd{expected  test error}
\kwd{random-$X$ regression}
\end{keyword}

\end{frontmatter}

\section{Introduction}
\label{sec:Introduction}

Regression is one of the most popular and well-known
methods in statistics, and it has been widely used across various fields.
Classical regression focuses on modeling the mean of a response variable $Y$ as a function of a set of explanatory variables $X$.  One common assumption
is to treat $X$  as fixed, i.e., nonrandom, which  \citet{Breiman1992}
and \citet{Rosset2020} refer to as ``Fixed-$X$''.  This assumption is reasonable if $X$ is constructed, such as in experimental design.
 Otherwise, due to real-world data collection and
predictive inferences, it may be more realistic to treat $X$ as random, namely ``Random-$X$''.
The seminal work of \citet{Rosset2020} provided a  detailed review of the literature for both the Fixed-$X$ and Random-$X$ settings in classical mean regression models. In addition, they analyzed the difference between expected  test errors under these two settings, and proposed   the estimators of  expected test errors for model assessment of mean regression models.

Instead of  classical mean regression, this paper focuses on
 covariance regression models, which link the covariance of
the $p$-dimensional response vector $Y_i$ to  the explanatory variables $X_i$
for $i=1,\cdots, n$. Note that  $X_i$ can be a vector,  matrix,  or  tensor, depending on the context of the study, although our main focus here is on the matrix form.
In addition, we assume that $Y_i$ and
$X_i$ are  random variables defined on the probability space $(\Omega,\mF,\rmP)$, where $\Omega$ is the sample space, $\mF$ is the $\sigma$-algebra of the subsets of $\Omega$, and $\rmP$ is a probability measure on $(\Omega,\mF)$.
To ensure the covariance regression model is well-defined,
we next introduce  the  $\sigma$-algebra $\mF_{\mathcal{X}}=\sigma\langle X_1,\cdots,X_n\rangle\subset \mF$,
which is generated by $X_1,\cdots, X_n$.
Then we propose the covariance regression model given below:
\be\label{eq:model}
\Cov[Y_i|\mF_{\mX}]= C_\beta(X_i),\textrm{ for }i=1,\cdots, n,
\ee
where $C_\beta(\cdot)$ is a pre-specified parametric and measurable map
characterized by a $K$-dimensional parameter vector $\beta=(\beta^{(1)},\cdots,\beta^{(K)})^\top$ such that $C_\beta(X_i)\in \mathcal{S}_p^+\subset\mathbb{R}^{p\times p}$ for all $i=1,\cdots,n$, and $\mathcal{S}_p^+=\{C\in \mathbb{R}^{p\times p}: C \textrm{ is positive definite}\}$.
This covariance regression model encompasses various covariance matrix structures such as sphericity, compound symmetry, and banding; see  examples given below.

\subsection{Fixed-X and Random-X settings for covariance regression models with five examples}\label{sec:review}

We introduce five  examples of covariance regression models. The first  is the linear covariance regression
model  $\Cov[Y_i|\mF_{\mathcal{X}}]=C_\beta(X_i)=\sum_{k=1}^{K}\beta^{(k)}X_i^{(k)}$,
where   $X_i=(X_i^{(1)},\cdots,X_i^{(K)})$,
$\beta=(\beta^{(1)}, \cdots,\beta^{(K)})^\top$, and the $X_i^{(k)}$s are $p\times p$
network adjacency matrices or similarity matrices constructed from their associated features.
The detailed structure of this model can be found in  \citet{Zou:Lan:Wang:Tsai:2017, Zou:Lan:Li:Tsai2021}.

The second example is the network autoregressive model, which can be used to describe the relationship between response variables of different network nodes.
Let $X_i=(X_{i,jj'})_{p\times p}$ be a  symmetric  $p\times p$ row-normalized network adjacency matrix with  $j,j'\in\{1,\cdots,p\}$
representing the $j$-th and $j'$-th nodes of the underlying network, and $X_{i,jj}=0$.
Then, the network autoregressive model has the  expression of $Y_{ij}=\beta^{(1)}\sum_{j'\neq j}X_{i,jj'} Y_{ij'}+\varepsilon_{ij}$ for $j=1,\cdots,p$, whose
vector-matrix form is $Y_i=\beta^{(1)}X_i Y_i+\mE_i$ ($i=1,\cdots,n$),  where the $Y_i=(Y_{i1},\cdots,Y_{ip})^\top$ consists of $p$ measurements associated with the $p$ nodes in the network,
 $\beta^{(1)}$ is the network autoregressive coefficient,
 $\mE_i=(\varepsilon_{i1},\cdots,\varepsilon_{ip})^\top$ is a $p$-dimensional random vector that is independent of  $\{X_i:i=1,\cdots,n\}$ with mean zero and covariance matrix $\sigma^2 I_p$, and $I_p$ is the $p$-dimensional identity matrix;
see, e.g., \citet{Zhou2017},  \citet{Huang2019} and \citet{Ma2019}. Based on this model, we have
 $\Cov[Y_i|\mF_{\mathcal{X}}]= C_\beta(X_i)=\sigma^2 (I_p-\beta^{(1)}X_i)^{-2}$, where $\beta=(\beta^{(1)}, \sigma^2)^\top$.

 The third example is the covariance regression model $\Cov[Y_i|\mF_{\mathcal{X}}]= C_\beta(X_i)=\Psi+BX_iX_i^\top B^\top$ proposed by \citet{Hoff2012}. For each $i=1,\cdots,n$,   $X_i$ is a $d$-dimensional random
covariate vector. In addition, $\Psi$ is a $p\times p$ nonrandom positive definite matrix, $B$ is a $p\times d$ nonrandom coefficient matrix, and 
$\beta$ consists of all the parameters in $\Psi$ and $B$. Based on \citet{Hoff2012}, this model has an interpretation as a type of random-effects model $Y_i=B X_i \gamma_i+\mE_i$,  where $\gamma_i$ is a one-dimensional  normally distributed  random coefficient with mean zero and variance 1, $\mE_i$ is a $p$-dimensional normally distributed random error vector 
with mean zero and covariance matrix $\Psi$, and $\{\gamma_i:i=1,\cdots, n\}$, $\{\mE_i:i=1,\cdots, n\}$ and $\{X_i:i=1,\cdots, n\}$ are mutually independent.

 The fourth example is that $C_\beta(X_i)$ is an  exponential function
of the random matrix (see \citealp{Chiu1996}). Finally, the fifth example is that
$C_\beta(X_i)$ is a function obtained from a Bayesian nonparametric covariance regression model (\citealp{Fox2015});  see Section \ref{sec:CheckingA1} in the supplementary material for details.

 The above five examples can be either Fixed-$X$ or Random-$X$ covariance 
models, depending on the settings of the $X_i$s. For example, the second example is a random model based on the $X$-setting in its cited literature. On the other hand, it can represent a spatial autoregressive model (\citealp{Lee2004} and \citealp{LeSage2009}), where $X_i=x_i$ is a spatial-weights matrix constructed by non-random geographical coordinates. Specifically,  $X_1=x_1,\cdots,X_n=x_n$ are fixed in this case, and $\mF_\mX$ and the covariance model reduce to the trivial $\sigma$-algebra  $\{\emptyset,\Omega\}$ and
$\Cov[Y_i|\mF_\mX]=\Cov(Y_i)= C_\beta(x_i)$, respectively. Similar examples include   $x_i$ having different structures of  geographical coordinates   in spatial statistics (see, e.g., \citealp{Cressie1993} and \citealp{Wikle2019}), and $x_i$ having a form of basis functions (see, e.g., \citealp{Scott2001}).

A special case of the Fixed-$X$ setting is where the $x_i$s are all equal to $x$
 for $i=1,\cdots, n$, which leads to $\Cov(Y_i)= C_\beta(x)$. This comprises various structured covariance matrices as  its special cases. For example, set $x=1_p1_p^\top-I_p$, where $1_p$ is the $p$-dimensional vector of ones. Then it yields the compound symmetry covariance matrix $\Cov(Y_i)= C_\beta(x)=\beta^{(1)}I_p+\beta^{(2)} x$, where $\beta=(\beta^{(1)}, \beta^{(2)})^\top$. Analogously, the banded covariance matrix can  be represented by the form $C_\beta(x)$ with a properly designed and fixed $x$ (see \citealp[p. 268]{Zou:Lan:Wang:Tsai:2017}).  
This special setting of $x$ constitutes another scenario of the  Fixed-$X$ setting.

 In contrast to the Fixed-$X$ setting, 
it is common to consider the $X_i$s as random along with  $Y_i$. 
Based on the literature in the above five examples, it is not surprising to find that these examples are Random-$X$ covariance models. For the sake of simplicity,  however, some of them assumed that, once observed, $X_i=x_i$ is fixed, and the ensuing theoretical results relied on that assumption.
The main difference between the fixed $X_i=x_i$ and the random $X_i$ settings is that, the fixed $X_i=x_i$ are 
 implicitly assumed to be independent for $i=1,\cdots,n$. This is because any nonrandom constants $x_{i}$ are
 mutually  independent for $i=1,\cdots,n$, which largely simplifies the derivation of asymptotic results for parameter estimators of covariance regression models. On the other hand,
$X_i$ is commonly dependent in real applications
 of  these examples.
Accordingly, this motivates us to propose the Random-$X$ covariance regression model, which can largely broaden the applications of these examples.

Based  on our understanding, there is no rigorous study of estimation theory under the Random-$X$ covariance setting. One possible reason  is that
the covariance regression model (\ref{eq:model}) under the Random-$X$  setting
induces an intrinsically  complicated dependence structure among the $p$ components of  $Y_i$. This can be even more challenging when  the $X_i$s are dependent for $i=1,\cdots,n$,  which commonly appears in the  literature; see, e.g.,  the above five examples. Consequently, it becomes a challenging task to develop  the asymptotic properties of  parameter estimators, in particular when $p\to\infty$.

As for the impact of Fixed-$X$ and Random-$X$ on the model assessment theory and expected test errors, it has only been studied
 in a seminal paper of  \citet{Rosset2020},  who considered a univariate response variable $Y_i$ with $p=1$ for mean regression models. In contrast, the covariance regression involves a multivariate response vector $Y_i$ with dimension $p>1$. Hence, it is intrinsically  challenging to develop  model assessment theory
based on covariance regression models. The usefulness for establishing this theory can be illustrated via  an empirical application in Section \ref{sec:empirical} of the supplementary material. This application is based on the linear covariance regression model; see the first example above.  By the assessment theory of covariance models,  one can select better models and  then make  accurate inferences. In Section \ref{sec:empirical}, we demonstrate that the selected model based on the developed  model assessment theory can  yield an optimal portfolio,
via the Markowitz optimization approach \citep{Markowitz1952},  which is  superior to the market portfolio  in terms of the Sharpe ratio. 
Some references on portfolio estimation or selection can be found in \citet{DeMiguel2009}, \citet{Fan2012}, and \citet{Cai2020}.

\subsection{Our contributions}
\label{sec:contribution}
The first aim of this paper is to construct a  framework for studying the estimation theory  under
the Random-$X$ covariance setting without imposing a normality assumption on the data while allowing the dimension $K$ to go to infinity.  Within this framework, we propose two types of parameter estimators for random covariance regression  (\ref{eq:model}). One is the quasi-maximum likelihood estimator (QMLE) and the other is the weighted least squares (WLS) estimator.
We then  demonstrate the consistency and asymptotic normality
of these two estimators under the following three scenarios: (i) $p$ is fixed and $n\to\infty$; (ii) $n$ is fixed and $p\to\infty$; and (iii) $n\to\infty$ and $p\to\infty$;  where $p$ is the dimension of the response vector $Y_i$, and $i=1,\cdots,n$.
Subsequently, we show that the resulting covariance matrix estimators are  consistent. We further compare the  asymptotic variances of the WLS estimator and QMLE, and demonstrate that the WLS estimator with the optimal weight matrix  is asymptotically the most efficient among all WLS estimators, and it is asymptotically as efficient as the QMLE under the normality assumption.

In sum, this paper has made three novel contributions to the estimation theory of Random-$X$ covariance regression.
(I) In deriving estimation theory, we do not impose any distributional  or independent and identically distributed (i.i.d.) conditions on $(Y_i,X_i)$, 
nor do we impose any independence condition on $X_i$, for $i=1,\cdots,n$. Instead, we only assume some moment conditions on $C_\beta(X_i)$ and its first three derivatives with respect to $\beta$,  and impose the  convergence in probability  conditions on  the averages of some functions of the random $X_i$s. The former conditions are  used to demonstrate the important distinctions between the Fixed-$X$ and Random-$X$ settings for estimation theory.
In addition, we allow a general error form for the data generating process, the model to be misspecified, and the conditional dependence of $Y_i$.
(II) The parametric form of
$C_\beta(X_i)$ in (\ref{eq:model}) can be linear or nonlinear, which contains various models in the literature as special cases, such as all the models introduced in
Section \ref{sec:review}. In addition, we allow  the dimensions of  $C_\beta(X_i)$ (i.e., $p$) and $\beta$ (i.e., $K$),  to go to infinity.  The formulations of (I) and (II)  are general under the Random-$X$ setting, and hence result in technical challenges and complexities in establishing estimation theory. 
(III) Compared to  classical large sample approaches (see, e.g., \citealp{White1982} and \citealp{van1998}) that commonly
assume that the $(Y_i,X_i)$ are  i.i.d. with fixed $p$  and $n\to\infty$,
 we need to develop a new  central limit theorem (CLT) and a new law of large numbers (LLN) for general
quadratic forms with random matrices to show our asymptotic results;  hence, new technical ingredients are constructed to demonstrate the new CLT and LLN.

In addition to establishing the estimation theory for the  Random-$X$ covariance regression model,
the second aim of this paper is to develop pioneering work on the 
model assessment theory. To this end,
we obtain the bias-variance decompositions for the expected test errors under both  Fixed-$X$ and Random-$X$  settings. This allows
us to propose model assessment criteria.
In sum, our contributions to  the model assessment theory are given below.
(I)  We  demonstrate that randomness in $X$ can yield additional bias and variance
in the  expected test errors
based on the WLS estimator. Furthermore,
treating the random $X$ as fixed can cause misleading estimators 
of expected test errors. (II) Developing model assessment theory of covariance regression with $p>1$ is more complicated than that of mean regression  with $p=1$ in \citet{Rosset2020};  taking into account  the complex dependence
structure from covariance regression results in technical challenges and complexities. Hence, we need to develop new approaches to demonstrate our model assessment  theories.

In conclusion, the pioneering work of this paper opens a new avenue in
studying random covariance regression models.
Based on the estimation and  model assessment theories developed under the Random-$X$ covariance setting, one can not only enhance the accuracy of estimation 
and model assessment, but can also explore new applications in the related fields.

\subsection{Organization of the paper}\label{sec:organization}
The remainder of this paper is organized as follows.
Section \ref{sec:2} establishes the estimation structure for the random covariance regression model.  Two parameter estimators, QMLE and WLS, are proposed and their corresponding theoretical properties are obtained.
Section \ref{sec:prediction} builds up the framework of model assessment, and analyzes the difference in expected test errors between the Random-$X$ and Fixed-$X$ settings. Section \ref{sec:40}  provides  simulation  studies and discussions. Section \ref{sec:5} gives concluding remarks. To save space, additional simulation results, an empirical study of stock returns in the US market, and additional tables and  all figures obtained from simulation studies and the empirical example, are presented in  Section \ref{sec:Ssimu}, Section \ref{sec:empirical}, and Section \ref{sec:tbfig} of the supplementary material, respectively.   All theoretical proofs, and some lemmas,  propositions, theorems,
assumptions and remarks  are also relegated to the supplementary material.

\section{Estimation and Inference with Random Covariance Regression}

\label{sec:2}

In this section, we propose two types of estimation approaches when the random covariance regression is correctly specified by the parametric model (\ref{eq:model}).
Denote the true covariance regression model  by
\be\label{eq:exp}
\Cov[Y_i|\mF_\mX]=C_{\beta_0}(X_i),\textrm{ for }i=1,\cdots, n,
\ee
where $\beta_0=(\beta_0^{(1)},\cdots,\beta_0^{(K)})^\top$ is the $K$-dimensional true parameter vector of $\beta$.
Since our main interest is covariance regression, the response $Y_i$ has been standardized to have conditional mean zero. Before studying $\beta$'s estimation and inference in the following subsections, we introduce some notation below. Let $\|\cdot\|_{\zeta}$ be the vector $\zeta$-norm or  the matrix $\zeta$-norm for $1\leq \zeta\leq  \infty$. Specifically, for any generic vector $\alpha=(\alpha_1,\cdots,\alpha_q)^\top\in\mathbb{R}^q$, $\|\alpha\|_{\zeta}=(\sum_{l=1}^q |\alpha_l|^{\zeta})^{1/{\zeta}}$, and, for any generic matrix $Z\in\mathbb{R}^{m\times q}$,
\be\label{eq:matrixnorm}
\laak Z\raak_{\zeta}=\sup\lbk \frac{\laak Z \alpha\raak_{\zeta}}{\laak \alpha\raak_{\zeta}}:\alpha\in\mathbb{R}^{q}\textrm{ and }\alpha\neq 0\rbk.
\ee
In addition, we utilize   the notation $Z=( z_{k_1k_2})_{m\times q}$ to denote an ${m}\times{q}$ matrix $Z$ with the $(k_1,k_2)$-th element
$z_{k_1k_2}$,
for  $k_1=1,\cdots, m$  and $k_2 = 1,\cdots,q$, and  let $|Z|=( |z_{k_1k_2}|)_{m\times q}$.   If $z_{k_1k_2}$ is a function of $\beta\in \mathfrak{B}$, namely $Z=Z(\beta)=\big( z_{k_1k_2}(\beta)\big)_{m\times q}$, then we denote $\sup_{\beta\in \mathfrak{B}}|Z(\beta)|=\big( \sup_{\beta\in \mathfrak{B}}|z_{k_1k_2}(\beta)|\big)_{m\times q}$.
If matrix $Z=(z_{k_1k_2})_{m\times q}$ is random, then, for $1\leq \zeta\leq\infty$,
denote
$\|z_{k_1k_2}\|_{L^\zeta}=(\rmE|z_{k_1k_2}|^{\zeta})^{1/{\zeta}}$
as the $L^{{\zeta}}$-norm of the  random variable $z_{k_1k_2}$.
Then define
\be\label{eq:Lsnorm}
\|Z\|_{L^\zeta}=\lsk \|z_{k_1k_2}\|_{L^\zeta}\rsk_{m\times q},
\ee
which is an $m\times q$ matrix  consisting of the $L^{{\zeta}}$-norm of
$z_{k_1k_2}$.
Moreover, denote
\be\label{eq:dnormmm}
\normmm{Z}_{L^{\zeta_1},\zeta_2}=\big\|\|Z\|_{L^{\zeta_1}}\big\|_{\zeta_2}
\ee
 for $\zeta_1,\zeta_2\in[1,\infty]$, where $\zeta_1$ and $\zeta_2$  correspond to the $L^{\zeta_1}$-norm (\ref{eq:Lsnorm}) and the matrix $\zeta_2$-norm (\ref{eq:matrixnorm}), respectively.
Lastly, let ${\rm diag}\{g_1,\cdots,g_p\}$ denote a generic diagonal matrix with diagonal elements  $g_1,\cdots, g_p$, and let $\circ$ represent the Hadamard product of any two matrices with the same dimensions.


\subsection{Quasi-maximum likelihood estimation under the Random-$X$  setting}

\label{sec:qmle}

Under parametric model (\ref{eq:model}),
the covariance coefficients $\beta$ can be estimated via the quasi-maximum likelihood estimation (QMLE) approach (e.g., see \citealp{Bollerslev1992} and \citealp{Tsay:2014}). We denote this estimator by $\hat\beta_{\QMLE}$, and it maximizes the log-likelihood function
\be\label{eq:likelihood}
\ell(\beta)=-\frac{np}{2}\log (2\pi)-\frac{1}{2}\sum_{i=1}^n{ \log\lmk {\rm det}\lbk  C_\beta(X_i)\rbk\rmk}-\frac{1}{2}\sum_{i=1}^n Y_{i}^\top  C^{-1}_\beta(X_i) Y_{i},
\ee
where ${\rm det}\{\cdot\}$ is the matrix determinant.
It is worth noting that even though this log-likelihood function is formulated
via the normal distribution, the consistency and asymptotic normality of QMLE $\hat \beta_{\QMLE}$ can be established without assuming that $Y_i$ is Gaussian distributed.

In order to establish the asymptotic normality of $\hat \beta_{\QMLE}$, we introduce the following assumption.


\begin{as}\label{as:1}
Assume that the data generating process is
$Y_{i}=C_{\beta_0}^{1/2}(X_i)\mE_{i}$ for
$i=1,\cdots,n,$
where  $C_{\beta_0}^{1/2}(X_i)$ is symmetric and the square root of covariance matrix $C_{\beta_0}(X_i)$, 
$\mE_{i}=(\varepsilon_{i1}, \cdots, \varepsilon_{ip})^\top,$
$\varepsilon_{ij}$ is
 conditionally independent  relative to the $\sigma$-algebra $\mF_{\mathcal{X}}$
 in both $i$ and $j$ (see the definition of conditional independence in \citealp[p. 322]{Chung2001}), $\rmE[\varepsilon_{ij}|\mF_{\mathcal{X}}]=0$,
and $\rmE[\varepsilon_{ij}^2|\mF_{\mX}]=1$ for all $i=1,\cdots,n$ and $j=1,\cdots,p$.
In addition, assume that  there exists a finite constant $\eta_1>0$ such that $\sup_{i,j}\rmE|\varepsilon_{ij}|^{4+\eta_1}<\infty$.
 \end{as}

Under Assumption \ref{as:1},  model (\ref{eq:exp}) is the true model and not misspecified. In addition, if the
$\mE_i$s are independent and also independent of $X_i$,
then $Y_{i}=C^{1/2}_{\beta_0}(X_i)\mE_{i}$ is similar to the independent component model (ICM, see, e.g., \citealp{Bai2004}, \citealp{Pan2008}, and \citealp{Zheng2012}) with  random and exogenous $X_i$ involved.  Although
the main results  presented in this section are based on the multiplicative error form $Y_{i}=C^{1/2}_{\beta_0}(X_i)\mE_{i}$, a more general error form for obtaining the estimation theory of model (\ref{eq:exp}) is given in Section \ref{sec:rAssumption1} of the supplementary material.  In addition, by Assumption \ref{as:1}, we  have that $\rmE[\mE_i\mE_i^\top|\mF_\mX]$ $=I_p$. Accordingly,
\be\label{eq:EYY}\rmE[Y_iY_i^\top|\mF_\mX]=C_{\beta_0}(X_i),
\ee which determines the  conditional dependence structure for the $p$ components in the response vector $Y_i$.  Furthermore,  Section \ref{sec:IntepretationA1} of the supplementary material demonstrates that $\Cov\big(\rmvec(Y_{i_1}Y_{i_1}^\top),\rmvec(Y_{i_2}Y_{i_2}^\top)\big)=\Cov( \rmvec \lbk C_{\beta_0}(X_{i_1})\rbk,\rmvec \lbk C_{\beta_0}(X_{i_2})\rbk)\neq 0$ for $i_1,i_2\in\{1,\cdots,n\}$ and $i_1\neq i_2$, where $\rmvec(G)$ denotes the vectorization for any generic matrix $G$; this is because we do not impose  any independent condition on $X_i$ for $i=1,\cdots,n$.   Accordingly, $Y_i$ is dependent for $i=1,\cdots, n$. Although  $Y_i$ is  conditionally independent  by Assumption \ref{as:1}, the estimation theory of model (\ref{eq:exp}) can be established by relaxing this conditional independence 
assumption; see Section \ref{sec:rconindpt} of the supplementary material.
More remarks on Assumption \ref{as:1},  such as  interpretations,  the examination of Assumption \ref{as:1} with regard to the five examples in Section \ref{sec:review}, and the significance and necessity of this assumption,    
are provided in   Section \ref{sec:RemarkA1} of the supplementary material.  We next introduce the assumption for $C_\beta(X_i)$ in model (\ref{eq:exp}).


 \begin{as} \label{as:2a} Assume that, for $i=1,\cdots,n$,
$C_\beta(X_i)$ denotes an indexed series of  positive definite  matrices, uniformly over $\beta$ that lies in a compact parameter
space  $\mathfrak{B}\subset\mathbb{R}^K$, $K<\infty$, and  $\beta_0$ is an interior point of  $\mathfrak{B}$. In addition, assume that 
$C_\beta(x_i)$ is three times differentiable with respect to $\beta$  for 
$P_{X_i}$-almost  every $x_i$ (see \citealp[p. 53]{van1998}), and
$C_\beta(x_i)$ and its first three derivatives with respect to $\beta$  are all measurable functions of $x_i$, where 
$P_{X_i}=\rmP\circ X_i^{-1}$ (see \citealp[p. 36]{Chung2001}) is the induced measure from the probability space $(\Omega,\mF,\rmP)$.
\end{as}

Based on Assumption \ref{as:2a}, the main results in this section focus on the dimension of  $\beta$, $K$, being fixed.
The estimation theory of model (\ref{eq:exp}) allowing  $K\to\infty$ is given in Section \ref{sec:KInf} of the supplementary material.
Under Assumptions \ref{as:1} -- \ref{as:2a},
one can derive the quasi-score function evaluated at $\beta_0$, which  is
\bea 
\frac{\partial\ell(\beta_0)}{\partial \beta^{(k)}}&=&\frac{1}{2}\sum_{i=1}^n Y_i^\top C^{-1}_{\beta_0}(X_i) \frac{\partial  C_{\beta_0}(X_i)}{\partial \beta^{(k)}} C^{-1}_{\beta_0}(X_i)  Y_i-\frac{1}{2}\sum_{i=1}^n\tr\lbk  C^{-1}_{\beta_0}(X_i)  \frac{\partial  C_{\beta_0}(X_i)}{\partial \beta^{(k)}}\rbk\nn\\
&=&\frac{1}{2}\sum_{i=1}^n \mE_i^\top C^{-1/2}_{\beta_0}(X_i) \frac{\partial  C_{\beta_0}(X_i)}{\partial \beta^{(k)}} C^{-1/2}_{\beta_0}(X_i)  \mE_i\label{eq:score11}-\frac{1}{2}\sum_{i=1}^n\tr\lbk  C^{-1}_{\beta_0}(X_i)  \frac{\partial  C_{\beta_0}(X_i)}{\partial \beta^{(k)}}\rbk
\eea
for $k=1,\cdots,K$, where, for any generic differentiable function $g(\beta)$ of $\beta$, we denote $\partial g(\beta)/\partial \beta^{(k)}$ evaluated at $\beta_0$ by $\partial g(\beta_0)/\partial \beta^{(k)}$ to save space. It is worth noting that Random-$X$ raises challenges for obtaining the asymptotic properties of the parameter estimators  under scenarios of  ``fixed $n$ and $p\to\infty$'' and ``$n\to\infty$ and $p\to\infty$''. For instance,
to derive the asymptotic normality of QMLE,
the central limit theorem (CLT)  needs to be developed  for $(np)^{-1/2}{\partial\ell(\beta_0)}/{\partial \beta^{(k)}}$.
Note that (\ref{eq:score11}) includes
a summation of the quadratic forms in $\mE_i$ with the random
matrices $C^{-1/2}_{\beta_0}(X_i)\{{\partial  C_{\beta_0}(X_i)}/{\partial \beta^{(k)}} \}C^{-1/2}_{\beta_0}(X_i)$, and $i=1,\cdots,n$.
Accordingly, the standard asymptotic theories for quadratic forms with nonstochastic matrices \citep{Kelejian2001} cannot be directly applied due to the random $X_i$.
In addition, to show the asymptotic normality of $\hat \beta_{\QMLE}$,  the law of large numbers (LLN) is needed for the
second derivative  evaluated at $\beta_0$,  i.e., $(np)^{-1} {\partial^2 \ell(\beta_0)}/(\partial \beta^{(k_1)}\partial \beta^{(k_2)})$, which also
contains the quadratic forms in $\mE_i$ with random matrices;
see equation (\ref{eq:partial2}) in the supplementary material.
To solve the problem, we propose a new CLT and a new LLN for
the summation of quadratic forms in $\mE_i$ with general random matrices.
Note that the  new CLT and LLN can be applied under any one of the following three scenarios: (i) $p$ is fixed and $n\to\infty$; (ii) $n$ is fixed and $p\to\infty$; and (iii) $n\to\infty$ and $p\to\infty$;
see Lemma \ref{cy:1} and its extensive proof  in Section \ref{sec:S2}  of
the supplementary material.   Furthermore,  the proof sketch for the CLT and  the new technical ingredients required in the proof are presented in Sections \ref{sec:LemmaS.10} of the supplementary material. By Lemma \ref{cy:1}, we demonstrate that $\hat\beta_{\QMLE}$ is consistent and asymptotically normal under the two additional assumptions  (i.e., Assumptions \ref{as:2b} -- \ref{as:3}) on the random matrices $C_\beta(X_i)$ given below. To save space,  the remarks of interpretations, examinations, significances, and necessities of Assumption \ref{as:2a}  as well as Assumptions \ref{as:2b} -- \ref{as:3} are given in  Sections \ref{sec:RemarkA2} -- \ref{sec:RemarkA4} of the supplementary material.


 \begin{as} \label{as:2b}Assume that
there exists finite $\eta_2>8/\eta_1$ such that
\bea
&&\max\biggl\{\sup_{p\geq 1}\sup_i\lrnormmm { C_{\beta_0}^{1/2}(X_i)}_{L^{\delta_c(2+\eta_2)},1},\sup_{p\geq 1}\sup_i\lrnormmm {\sup_{\beta\in\mathfrak{B}}\lak C_{\beta}^{-1}(X_i)\rak}_{L^{\delta_0(2+\eta_2)},1}, \label{eq:normmmc}\\
&&\sup_{p\geq 1}\sup_i\lrnormmm {\sup_{\beta\in\mathfrak{B}}\lak \frac{\partial C_\beta(X_i)}{\partial \beta^{(k)}}\rak}_{L^{\delta_1(2+\eta_2)},1},
\sup_{p\geq 1}\sup_i\lrnormmm {\sup_{\beta\in\mathfrak{B}}\lak \frac{\partial^2 C_\beta(X_i)}{\partial \beta^{(k_1)}\partial \beta^{(k_2)}}\rak}_{L^{\delta_2(2+\eta_2)},1},\nn\\
&&\sup_{p\geq 1}\sup_i\lrnormmm { \sup_{\beta\in\mathfrak{B}}\lak \frac{\partial^3 C_\beta(X_i)}{\partial \beta^{(k_1)}\partial \beta^{(k_2)}\partial \beta^{(k_3)}}\rak}_{L^{\delta_3(2+\eta_2)},1}\biggl\}<\infty,\nn
\eea
for any $k,k_1,k_2,k_3\in\{1,\cdots,K\}$, and for some $\delta_c,\delta_0,\delta_1,\delta_2,\delta_3\in (0,\infty)$ that satisfy
\[
 \max\lbk \lsk \frac{2}{\delta_c}+\frac{4}{\delta_0}+\frac{3}{\delta_1}\rsk,\lsk \frac{2}{\delta_c}+\frac{3}{\delta_0}+\frac{1}{\delta_1}+\frac{1}{\delta_2}\rsk,\lsk \frac{2}{\delta_c}+\frac{2}{\delta_0}+\frac{1}{\delta_3}\rsk\rbk\leq 1,
\]
where $\eta_1$ is defined in Assumption \ref{as:1}.
 \end{as}

Before introducing the next assumption,
 define \be\label{eq:aiq}
A_{i,{\rm Q}}^{(k)}=\frac{1}{2}C_{\beta_0}^{-1/2}(X_i)\frac{\partial C_{\beta_0}(X_i)}{\partial \beta^{(k)}}C_{\beta_0}^{-1/2}(X_i)\ee and $\Upsilon_i^{(4)}={\rm diag}\{\mu^{(4)}_{i  1}-3,\cdots, \mu^{(4)}_{i  p}-3\}\in\mathbb{R}^{p\times p}$, where $\mu^{(4)}_{i  j}=\rmE[\varepsilon_{ij}^4| \mF_{\mX} ]$,  and ``${\rm Q}$'' used in the subscript is an abbreviation of ``${\rm QMLE}$'' to save space. In addition, denote
$
V_{np,{\rm Q}}={2}({np})^{-1}\sum_{i=1}^n\big(\tr(A_{i,{\rm Q}}^{(k_1)}A_{i,{\rm Q}}^{(k_2)}) \big)_{K\times K}\textrm{ and }\Omega_{np,{\rm Q}}=({np})^{-1}\sum_{i=1}^n \big(\tr(A_{i,{\rm Q}}^{(k_1)}\circ \Upsilon_i^{(4)}\circ A_{i,{\rm Q}}^{(k_2)}) \big)_{K\times K}.
$


 \begin{as}\label{as:3}  As $\max\{n,p\}\to\infty$, we have that $V_{np,{\rm Q}}\stackrel{\rmP}\longrightarrow V_{\rm Q}$ and $\Omega_{np,{\rm Q}}\stackrel{\rmP}\longrightarrow \Omega_{\rm Q}$, where   $V_{\rm Q}$ and $\Omega_{\rm Q}$ are finite $K\times K$ matrices. In addition, assume both
  $V_{\rm Q}$ and $V_{\rm Q}+\Omega_{\rm Q}$  are positive definite.
 \end{as}

It is worth noting that Assumptions \ref{as:1} -- \ref{as:3} do not impose any distributional  and  i.i.d. conditions on $(Y_i,X_i)$, or an independent condition on $X_i$, for $i=1,\cdots,n$. Instead, we only assume some moment conditions on $C_\beta(X_i)$ and its first three derivatives with respect to $\beta$ in   Assumption \ref{as:2b}, and impose the convergence in probability  conditions on  the averages of some functions of random $X_i$s in Assumption \ref{as:3}. The conditions in Assumption \ref{as:2b} can be used to demonstrate the important distinctions between the Fixed-$X$ and Random-$X$ settings for estimation theory; see  Remark \ref{eq:advantage} in  Section \ref{sec:CheckingA3} of the supplementary material.  
In  the two examples   demonstrating the justification for Assumptions \ref{as:2b} -- \ref{as:3}  presented in Section \ref{sec:analytical} of the supplementary material,  we find that, although Assumption \ref{as:3} does not directly involve a distributional condition on $X_i$, the justification for Assumption \ref{as:3} possibly requires some distributional and finite moment conditions on the functions of $X_i$. 
More remarks on Assumptions  \ref{as:2b} and  \ref{as:3}  are presented in  Sections  \ref{sec:RemarkA3} and  \ref{sec:RemarkA4}, respectively,  of the supplementary material.

The above assumptions lead to the following conclusion.


 \begin{tm} \label{tm:2}Under Assumptions \ref{as:1} -- \ref{as:3},  we have that, as $\max\{n,p\}\to\infty$,
\[
\sqrt{np}\lsk \hat\beta_{\QMLE}-\beta_0\rsk\stackrel{d}\longrightarrow N(0, V_{\rm Q}^{-1}+ V_{\rm Q}^{-1} \Omega_{\rm Q} V_{\rm Q}^{-1}).
\]
\end{tm}

\noindent In Theorem \ref{tm:2} and the theorems developed below, we unify the following three scenarios by $\max\{n,p\}\to\infty$: (i) $p$ is fixed and $n\to\infty$; (ii) $n$ is fixed and $p\to\infty$; and (iii) $n\to\infty$ and $p\to\infty$.
 Note that the proofs of all the theorems and propositions in this paper are presented in Sections \ref{sec:pnproofs} -- \ref{sec:thmproofs} of the supplementary material. In order to help better understand the key steps of showing the theorems and propositions, we have provided the flow charts and discussions in  Section \ref{sec:flowchart0} of the supplementary material.
In addition to the QMLE, we next consider a least squares type estimator.

\subsection{Weighted least squares estimation under the Random-$X$  setting}
\label{sec:lcr}

Without the distribution assumption, it is natural to
consider the weighted least squares (WLS) approach to obtain
parameter estimators.
To this end,  let  $\bmY_i=\rmvec(Y_iY_i^\top)$ and
$\bmC_\beta (X_i)=\rmvec\, C_\beta(X_i)$.
For the given $\beta$,
we  define  a $p^2$-dimensional error vector $\bmY_i-\bmC_\beta (X_i)$.  In order to aggregate each component of $\bmY_i-\bmC_\beta (X_i)$ and obtain a total error,  we further introduce a $p^2\times p^2$ positive semidefinite weight matrix  $\bmW(X_i)$, which depends on $X_i$. Afterwards, we define the total error by
\be\label{eq:Tr}
Tr(\beta)=\frac{1}{n}\sum_{i=1}^n \lbk \bmY_i- \bmC_\beta(X_i)\rbk^\top  \bmW(X_i) \lbk \bmY_i- \bmC_\beta(X_i)\rbk.
\ee
Accordingly, the WLS estimator $\hat\beta_{\WLS}$ can be obtained by
minimizing $Tr(\beta)$.
First let $\bmW(X_i)=I_{p^2}$. Then minimizing  (\ref{eq:Tr}) leads to the ordinary least squares (OLS) estimator. To improve the efficiency of the OLS estimator,
 let $\bmW(X_i)=(\Cov[\bmY_i|X_i])^-$ be the generalized (Moore-Penrose) inverse of $\Cov[\bmY_i|X_i]$.
We then obtain a generalized least squares (GLS)  estimator $\hat\beta_{\rm GLS}$  by minimizing (\ref{eq:Tr}). Theorem \ref{tm:asyeff} in Section \ref{sec:aeftge} of the supplementary material shows that $\hat\beta_{\rm GLS}$ is asymptotically the most efficient among all $\hat\beta_\WLS$.  It is worth noting that $(\Cov[\bmY_i|X_i])^-$ is only positive semidefinite but not positive definite, which motivates us to consider a positive semidefinite matrix $\bmW(X_i)$ in obtaining the WLS estimator.


 \begin{re}\label{re:m1}
 The closed form of $\Cov[\bmY_i|X_i]=\bmV(X_i)$ in finding the GLS estimator
 is given in Theorem \ref{tm:ccov} of Section \ref{subsec:SC} of the supplementary material, and it depends on $\beta_0$. Substituting $\beta_0$ with the OLS estimator
yields  an estimator  $\hat \bmV(X_i)$  of $\bmV(X_i)$; the details  of obtaining $\hat \bmV(X_i)$
are given in Section \ref{subsec:SC} of the supplementary material.
Following the GLS standard procedure  in  classical mean regression
 models (see, e.g., \citealp{Wooldridge2012}),   we then obtain a feasible GLS (FGLS) estimator by
minimizing  (\ref{eq:Tr}) with $\bmW(X_i)=\hat \bmV^-(X_i)$.
In fact, this procedure can be applied for any $\bmW(X_i)$ that depends on $\beta_0$.
\end{re}

Analogous to the derivation of the asymptotic normality of QMLE, the asymptotic property of WLS relies on the first- and second-order derivatives of $Tr(\beta)$  evaluated at $\beta_0$, i.e.,  $\partial \,Tr(\beta_0)/\partial \beta$ and $\partial^2 Tr(\beta_0)/(\partial \beta\partial \beta^\top)$. Since both  derivatives involve a summation of  quadratic forms with random matrices, Lemma \ref{cy:1}  in Section \ref{sec:S2} of the supplementary material is used to establish the asymptotic normality of $\hat \beta_{\WLS}$, which is given below. Note that the required Assumptions \ref{as:2b}$'$ and \ref{as:3}$'$ for WLS are parallel to
 Assumptions \ref{as:2b} and \ref{as:3} for QMLE. To save space, they are presented in   Section \ref{sec:A1additional} of the supplementary material.


\begin{tm} \label{tm:3}Under  Assumptions \ref{as:1}, \ref{as:2a},  \ref{as:2b}$\,'$ and \ref{as:3}$\,'$, we have that, as $\max\{n,p\}\to\infty$,
\[
\sqrt{np}\lsk \hat\beta_{\WLS}-\beta_0\rsk\stackrel{d}\longrightarrow N \big(0, M_{\rmG}^{-1}(V_{\rmG}+\Omega_{\rmG} )M_{\rmG}^{-1}\big),
\]
where   $V_{\rmG}$, $\Omega_{\rmG}$ and $M_{\rmG}$ are $K\times K$ matrices defined in Assumption \ref{as:3}$\,'$, and ``${\rm W}$'' used in the subscript is an abbreviation of ``${\rm WLS}$'' to save space.
\end{tm}

It is worth noting that Theorem \ref{tm:3} is based on a known weight matrix $\bmW(X_i)$, and we can demonstrate that this theorem is also valid for the FGLS estimator mentioned in  Remark \ref{re:m1}.


 \begin{re}\label{re:comparison}  Theorems \ref{tm:2} and \ref{tm:3} can lead to the asymptotic variances of $\hat\beta_{\QMLE}$, $\hat\beta_{\WLS}$, and $\hat\beta_{\GLS}$,
and they are, correspondingly,  $\Avar(\hat\beta_{\QMLE})=V_{\rm Q}^{-1}+ V_{\rm Q}^{-1} \Omega_{\rm Q} V_{\rm Q}^{-1}$,
$
\Avar(\hat\beta_{\WLS})$ $= M_{\rmG}^{-1}(V_{\rmG}+\Omega_{\rmG})M_{\rmG}^{-1},
$
and
$
\Avar(\hat\beta_{\GLS})= \{M_{\rmG}^{-1}(V_{\rmG}+\Omega_{\rmG})M_{\rmG}^{-1}\}{|}_{\bmW(X_i)=\bmV^-(X_i)}$.  
 We then compare these  asymptotic variances, and show that the GLS estimator  is asymptotically the most efficient among all WLS estimators, and it is asymptotically as efficient as the QMLE under the normality assumption; see Theorems \ref{tm:asyeff} -- \ref{tm:QMLEGLS} in Section \ref{sec:aeftge} of the supplementary material.  The connection between the  above optimality of 
GLS   and the standard    GLS' optimality result in classical mean regression models  (see, e.g., \citealp{Wooldridge2012}) is given in Remark \ref{re:link} of the supplementary material.
\end{re}


 \begin{re} \label{re:positived} Based on the above estimation approaches, the covariance matrix estimator can be obtained by substituting $\beta$ in $C_\beta(X_i)$ with $\hat\beta$, where $\hat\beta$ can be either $\hat\beta_{\QMLE}$ or $\hat\beta_{\WLS}$. However,
 to assure the positive definiteness of the covariance matrix estimator, we need to impose the constraint $\beta\in B^+=\{ \beta:C_\beta(X_i)\succ 0,\textrm{ for all }i=1,\cdots, n\}$, where   $C_\beta(X_i)\succ 0$ indicates that $C_\beta(X_i)$ is positive definite. Under this constraint, we obtain the constrained  estimator $\hat \beta^c=\arg\min_{\beta\in B^+}M(\beta)$, where $M(\beta)$ can be $-\ell(\beta)$ and $n\,Tr(\beta)$ used for obtaining the  constrained QMLE and WLS estimators, respectively.
The algorithm for obtaining the constrained estimator is presented in Section \ref{sec:algorithm} of the supplementary material.
In addition, employing Theorems \ref{tm:2} -- \ref{tm:3}, we can show that the
unconstrained estimator is identical to the constrained estimator with probability tending to one, as $\max\{n,p\}\to\infty$; see the proof after proving Theorem \ref{tm:3} in Section \ref{sec:thmproofs} of the supplementary material. Consequently, $\hat\beta$ and $\hat \beta^c$ are asymptotically equivalent.
 \end{re}


 \begin{re} \label{re:mis} Theorems \ref{tm:2} -- \ref{tm:3} are established under Assumption \ref{as:1}, which indicates that  model (\ref{eq:exp}) is the true model and not misspecified. When model (\ref{eq:exp}) is misspecified, we have demonstrated that the proposed QMLE and WLS estimators are asymptotically normal; see  Theorems \ref{tm:2M} -- \ref{tm:3M} in Section \ref{sec:misspecification} of the supplementary material. The numerical studies of parameter estimation under the misspecified model are presented in Section \ref{sec:Ssimu} of the supplementary material.
 \end{re}

Recall that the first example in  Section \ref{sec:review} introduces a linear covariance regression model,
\be\label{eq:lcov}
C_{\beta}(X_i)=\sum_{k=1}^K\beta^{(k)}X_i^{(k)}.
\ee
Under this model,
both of  \citeauthor{Zou:Lan:Wang:Tsai:2017}'s (\citeyear{Zou:Lan:Wang:Tsai:2017})
estimators (OLS and FGLS)
are special cases of our proposed WLS estimator.
Denote $\bmX_i=\big(\rmvec(X_i^{(1)}),\cdots,$ $\rmvec(X_i^{(K)})\big)$. Then, the resulting WLS estimator has a closed form
\be\label{eq:GLS}
\hat\beta_{\WLS}=\lbk \sum_{i=1}^n \bmX_i^\top \bmW(X_i) \bmX_i\rbk^{-1}\sum_{i=1}^n \bmX_i^\top \bmW(X_i) \bmY_i.
\ee
The WLS reduces to OLS by setting $\bmW(X_i)=I_{p^2}$, and the resulting OLS estimator has a closed  form
\[
\hat\beta_{\OLS}=\lsk \sum_{i=1}^n \bmX_i^\top\bmX_i\rsk^{-1}\sum_{i=1}^n \bmX_i^\top\bmY_i=\lsk\sum_{i=1}^n\textrm{tr}( X^{(k_1)}_iX^{(k_2)}_i ) \rsk_{K\times K}^{-1}\lsk \sum_{i=1}^nY_i^\top X^{(k)}_i Y_i\rsk_{K\times 1}.
\]

Based on  the theoretical properties of $\hat\beta_{\QMLE}$ and $\hat\beta_{\WLS}$,
 we next show the stochastic convergence rate of the
estimated covariance matrix.

\subsection{Inference with Random-$X$ covariance regression}

Applying the results of Theorems \ref{tm:2} and \ref{tm:3}, we can make inferences on the conditional covariance matrix
$\Cov[Y_i|\mF_\mX]=C_{\beta_0}(X_i)$.
Let $\|G\|_F$ be the Frobenius norm for any generic matrix $G$. To study the asymptotic properties of  the estimators of $C_{\beta_0}(X_i)$
and its inverse under the Frobenius norm,  we introduce an assumption below.


\setcounter{as}{4}

  \begin{as}\label{as:4}  For any  fixed $i$, as $p\to\infty$, assume that
  \[  p^{-1}\lsk {\rm tr}\lbk C_{\beta_0}^l(X_i) \frac{\partial C_{\beta_0}(X_i)}{\partial \beta^{(k_1)}} C_{\beta_0}^l(X_i) \frac{\partial C_{\beta_0}(X_i)}{\partial \beta^{(k_2)}}\rbk\rsk_{K\times K}\stackrel{\rmP}\longrightarrow  Q_{l,i},
  \]
for $l=-2$ and $l=0$, where $C_{\beta_0}^0(X_i)\triangleq I_p$  and  $i=1,\cdots,n$.
\end{as}
\noindent This assumption
 is a generalization of Condition (C3) from \citet{Zou:Lan:Wang:Tsai:2017} to random matrices $ C_{\beta_0}(X_i)$ and ${\partial C_{\beta_0}(X_i)}/{\partial \beta^{(k)}}$,
 and it is analogous to the standard condition  imposed on the asymptotic covariance matrix of the OLS  estimator obtained from classical mean regression models under the Random-$X$ setting.

We next present the theoretical properties of the estimators of $C_{\beta_0}(X_i)$
and the precision matrix given below, under Assumptions \ref{as:1} -- \ref{as:2b} and Assumption \ref{as:4}.


\begin{tm}  \label{tm:region}  Under Assumptions \ref{as:1} -- \ref{as:2b} and Assumption \ref{as:4}, for any fixed $i$, as $p\to\infty$, and with $n$ either  fixed or going to infinity, we have that
\[
n\laak C_{\hat\beta}(X_i)- C_{\beta_0}(X_i)\raak_F^2\stackrel{d}\longrightarrow Z^\top
Q_{0,i} Z,\textrm{ and }n\laak C_{\hat\beta}^{-1}(X_i)- C_{\beta_0}^{-1}(X_i)\raak_F^2\stackrel{d}
\longrightarrow Z^\top Q_{-2,i} Z,
\]
 given that $C_{\hat\beta}(X_i)$ is invertible a.s.,
where $\hat\beta$ is an estimator satisfying $\sqrt{np}( \hat\beta-\beta_0)\stackrel{d}\longrightarrow Z$ with $Z$ following a $K$-dimensional multivariate normal distribution.
\end{tm}

Theorem \ref{tm:region} indicates that the orders of  $p^{-1/2}\|C_{\hat\beta}(X_i)- C_{\beta_0}(X_i)\|_F$ and $p^{-1/2}\|C^{-1}_{\hat\beta}(X_i)- C^{-1}_{\beta_0}(X_i)\|_F$ are of $(np)^{-1/2}$.
Under the spectral norm $\|\cdot\|_2$, we  can  obtain the orders of $\|C_{\hat\beta}(X_i)- C_{\beta_0}(X_i)\|_2$ and $\|C^{-1}_{\hat\beta}(X_i)- C^{-1}_{\beta_0}(X_i)\|_2$ given below.


\begin{tm}  \label{tm:cons}  Under Assumptions \ref{as:1} -- \ref{as:2b}, for any fixed $i$,  we have  that,  as $\max\{n,p\}\to\infty$,
\[
\laak C_{\hat\beta}(X_i)- C_{\beta_0}(X_i)\raak_2=O_{\rmP}\lsk\frac{1}{\sqrt{np}} \rsk,\textrm{ and }\laak C^{-1}_{\hat\beta}(X_i)- C^{-1}_{\beta_0}(X_i)\raak_2=O_{\rmP}\lsk\frac{1}{\sqrt{np}}\rsk,
\]
given that $C_{\hat\beta}(X_i)$ is invertible a.s.,  where $\hat\beta$ is described in Theorem \ref{tm:region}.
\end{tm}

\noindent
Based on Theorems \ref{tm:2} and \ref{tm:3}, the estimators $\hat\beta_{\QMLE}$ and $\hat\beta_{\WLS}$ satisfy  the asymptotic distribution assumptions in Theorems \ref{tm:region} and \ref{tm:cons}, respectively. Thus, the theoretical results in Theorems \ref{tm:region} and \ref{tm:cons} are respectively applicable to these estimators.
In addition, Theorems \ref{tm:2} -- \ref{tm:cons} are developed under the Random-$X$ setting. Since Fixed-$X$ is a special case of Random-$X$, the asymptotic properties in Theorems \ref{tm:2} -- \ref{tm:cons} still hold under the Fixed-$X$ setting.

It is worth noting that the asymptotic variances of  the
QMLE and the WLS
estimator in Theorems \ref{tm:2} and \ref{tm:3}, respectively, are unknown in practice.
Hence, we need to find their corresponding estimators.
 For the sake of convenience, let $\hat\beta$ denote the consistent  estimator of $\beta_0$. In both Theorems \ref{tm:2} and \ref{tm:3},
we also
let $V$, $\Omega$ and $M$ denote the quantities that appeared in the asymptotic variance,
$\Avar (\hat\beta)$.
To estimate $V$, $\Omega$ and $M$,  we replace $\beta_0$  in (\ref{eq:aiq}) and (\ref{eq:aig})  in  Section \ref{sec:A1additional} of the supplementary material by  $\hat\beta$, which yields the corresponding estimators $\widehat A_{i,{\rm Q}}^{(k)}$ and  $\widehat A_{i,{\rmG}}^{(k)}$. Moreover, let
 $\widehat \Upsilon_i^{(4)}={\rm diag}\{\hat\varepsilon_{i1}^4-3,\cdots, \hat\varepsilon_{ip}^4-3\}\in\mathbb{R}^{p\times p}$,
where $(\hat \varepsilon_{i1},\cdots, \hat \varepsilon_{ip})^\top\triangleq\widehat \mE_{i}
=C^{-1/2}_{\hat\beta}(X_i) Y_{i}$.
Then, we obtain the estimators of
 $V$, $\Omega$ and $M$,  which are given below:
\bea\label{eq:hatvo}&&\\
\widehat V_{ np,{\rm Q}}&=&\frac{2}{np}\sum_{i=1}^n\lsk\tr(\widehat A_{i,{\rm Q}}^{(k_1)}\widehat A_{i,{\rm Q}}^{(k_2)}) \rsk_{K\times K},
\widehat\Omega_{ np,{\rm Q}}=\frac{1}{np}\sum_{i=1}^n \lsk\tr(\widehat A_{i,{\rm Q}}^{(k_1)}\circ \widehat \Upsilon_i^{(4)}\circ \widehat A_{i,{\rm Q}}^{(k_2)}) \rsk_{K\times K},\nn\\
\widehat V_{ np,{\rmG}}&=&\frac{2}{np}\sum_{i=1}^n\lsk\tr(\widehat A_{i,{\rmG}}^{(k_1)}\widehat A_{i,{\rmG}}^{(k_2)}) \rsk_{K\times K},
\widehat\Omega_{ np,{\rmG}}=\frac{1}{np}\sum_{i=1}^n \lsk\tr(\widehat A_{i,{\rmG}}^{(k_1)}\circ \widehat \Upsilon_i^{(4)}\circ \widehat A_{i,{\rmG}}^{(k_2)}) \rsk_{K\times K},\nn
\eea
\[
\textrm{and }\widehat M_{np,{\rmG}}=\frac{1}{np}\sum_{i=1}^n \lsk \rmvec^\top\lbk \frac{\partial C_{\hat\beta}(X_i)}{\partial \beta^{(k_1)}}\rbk  \bmW(X_i) \rmvec\lbk \frac{\partial C_{\hat\beta}(X_i)}{\partial \beta^{(k_2)}}\rbk\rsk_{K\times K}.\nn
\]

Note that Theorem \ref{tm:region} and Theorem \ref{tm:cons} separately show that
$C_{\hat\beta}(X_i)$ and $C^{-1}_{\hat\beta}(X_i)$ are consistent estimators of $C_{\beta_0}(X_i)$ and $C^{-1}_{\beta_0}(X_i)$.
Using these results, we can verify that
$
\widehat\Avar (\hat\beta_{\QMLE})= \widehat V_{ np,{\rm Q}}^{-1}+ \widehat V_{ np,{\rm Q}}^{-1} \widehat\Omega_{ np,{\rm Q}} \widehat V_{np,{\rm Q}}^{-1}$ and $\widehat\Avar$ $ (\hat\beta_{\WLS})= \widehat M_{np,\rmG}^{-1}\widehat V_{np,\rmG}\widehat M_{np,\rmG}^{-1}+ \widehat M_{np,\rmG}^{-1} \widehat \Omega_{np,\rmG} \widehat M_{np,\rmG}^{-1}
$
are consistent estimators of $\Avar (\hat\beta_{\QMLE})$ and $\Avar (\hat\beta_{\WLS})$,  respectively, which are the asymptotic variances of $\hat\beta_{\QMLE}$ and
$\hat\beta_{\WLS}$.
Using them, we can make inferences regarding $\beta$ such as constructing confidence intervals and testing hypotheses.

Finally,  the quantities $Q_{0,i}$ and $Q_{-2,i}$, which appear in the asymptotic distributions of Theorem \ref{tm:region}, are unknown. In practice, they can be replaced by the consistent estimators.
 Specifically, replacing $\beta_0$, which appears in Assumption \ref{as:4}, by $\hat\beta_{\QMLE}$ (or $\hat\beta_{\WLS}$),
yields the consistent estimators of $Q_{0,i}$ and $Q_{-2,i}$.
After studying parameter estimation and  inference, we next focus 
on  model  assessment.



\section{Model Assessment for Random-$X$  Covariance Regression}

\label{sec:prediction}

To assess the performance of Random-$X$ covariance regression models, we adopt the model assessment approach described in Chapter 7 of \citet{Hastie2009} to analyze  the training error and test error. Based on the training data 
$(\mbX,\mbY)=\big\{(X_i,Y_{i}):i \in\{1,\cdots,n\}\big\}$,
we can  estimate the unknown covariance regression parameter $\beta$ by $\hat\beta$, where  $\hat\beta$  can be either  $\hat\beta_{\QMLE}$ or $\hat\beta_{\WLS}$ proposed in  Section \ref{sec:2}. Recall
that in Section \ref{sec:lcr}, we denoted   $\bmC_\beta (X_i)=\rmvec\, C_\beta(X_i)$ and $\bmY_i=\rmvec(Y_i Y_i^\top)$ for the training data $i\in\{1,\cdots,n\}$.
 For the sake of simplicity, we let $\hat \bmC(\cdot)=\bmC_{\hat\beta}(\cdot)$,
 while  $\hat \bmC(\cdot)$ is sometimes denoted by $\hat \bmC_{\mbX,\mbY}(\cdot)$
to avoid confusion.  
To assess the model fit, we adopt (\ref{eq:Tr}) and define the training error
as follows:
 \be\label{eq:traininge}
\Tr(\mbX,\mbY)=\frac{1}{n}\sum_{i=1}^n L\big( \bmY_i, \hat \bmC(X_i)\big),
\ee
where the loss function $L\big( \bmY_i, \hat \bmC(X_i)\big)=\{ \bmY_i-\hat \bmC(X_i)\}^\top  \bmW(X_i) \{ \bmY_i-\hat \bmC(X_i)\}$  measures the error between $ \bmY_i$ and $\hat \bmC(X_i)$. Unfortunately, the training error $\Tr(\mbX,\mbY)$ cannot be used to assess the performance of covariance regression models; this is because the training error does not properly account for model 
complexity \citep[p. 38]{Hastie2009}. To this end, we propose the test error for covariance regression models, which is the average loss over the test data set $(\mbX_0,\mbY_0)=\big\{(X_{0i},Y_{0i}):i \in\{1,\cdots,m\}\big\}$, i.e.,
\be\label{eq:teste}
\Te(\mbX_0,\mbY_0,\mbX,\mbY)=\frac{1}{m}\sum_{i=1}^{m} L\big( \bmY_{0i}, \hat \bmC(X_{0i})\big),
\ee
where 
$\bmY_{0i}=\rmvec(Y_{0i} Y_{0i}^\top)$ and   
$Y_{0i}\in\mbY_0=\big\{Y_{0i}:i \in\{1,\cdots,m\}\big\}$.    
 It is worth noting that the aim of this section only focuses on  model assessment rather than utilizing $\hat \bmC(X_{0i})$ to predict
$\bmY_{0i}$. Hence, our focus is to assess the covariance regression models by analyzing the test error (\ref{eq:teste}) and its expected value given below. 
In addition, 
the test error $\Te(\cdot)$ depends on both  training data $(\mbX,\mbY)$ and  test data
$(\mbX_0,\mbY_0)$, where $(\mbX_0,\mbY_0)$ is not used to train the model $\hat \bmC(\cdot)=\bmC_{\hat\beta}(\cdot)$,   and we employ the
same
$\hat \bmC(\cdot)$ for both training and
test errors.  This is similar to the definitions of training and
test errors  in  classical mean regression models.

In  practice, test observations $(\mbX_0,\mbY_0)$ may not be available. In this case, we can only employ the training data $(\mbX,\mbY)$ to estimate the expected test error
\be\label{eq:prerr}
\Err=\rmE\lbk \Te(\mbX_0,\mbY_0,\mbX,\mbY)\rbk,
\ee
 and then use the estimator of $\Err$ to assess covariance regression models. This is similar to treating Mallows' Cp \citep{Mallows1973} and other model assessment criteria as  the estimators of the expected test error in  classical mean regression models (see, e.g., \citealp[p. 230]{Hastie2009}). 
In these models, \citet{Rosset2020} found that 
the Random-$X$ and Fixed-$X$ settings can result in different expected test errors. This motivates us to analytically evaluate such differences
for covariance regression models. To this end, we first study the  expected  test errors
under both  Random-$X$ and Fixed-$X$ settings. Then,  we make comparisons between these two  errors.  We next adapt the definitions of
 \citet{Rosset2020} and consider the following three settings for both 
training data $\mbX$ and test data $\mbX_0$.

\textbf{Fixed-$X$ Setting.} Assume that the explanatory variables in the training data
 set are fixed, i.e., $X_1=x_1,\cdots,X_n=x_n$.
In addition,  assume that the explanatory variables in the test data are the same as those in the
training data, i.e., $X_{01}=x_{1},\cdots,X_{0n}=x_{n}$.

\textbf{Random-$X$ Setting.} Assume that $(X_1,Y_1),\cdots,(X_n,Y_n), (X_{01},Y_{01}),\cdots,$ $(X_{0m},Y_{0m})$ are  i.i.d. according to some joint distribution of $(X,Y)$.

\textbf{Same-$X$ Setting.} Assume that $(X_1,Y_1),\cdots,(X_n,Y_n)$ are
 i.i.d. according to some joint distribution of $(X,Y)$.
In addition, assume that $X_{0i}=X_i$ for $i=1,\cdots,m$ with $m=n$, and the corresponding responses $Y_{0i}$ are independent distributed according to the conditional distribution of $Y_i|X_i$ for $i=1,\cdots,n$.

It is clear that both Random-$X$ and Same-$X$ settings have the random training data $\mbX=\big\{X_i:i \in\{1,\cdots,n\}\big\}$
and the random test data $\mbX_0=\big\{X_{0i}:i \in\{1,\cdots,m\}\big\}$.
Without causing confusion, we
also use $\big\{(X_\alpha,Y_{\alpha}):\alpha \in\{1,\cdots,n\}\big\}$ and $\big\{(X_\alpha,Y_{\alpha}):\alpha \in\{01,\cdots,0m\}\big\}$ to denote the training data and the test data, respectively.

We next evaluate the expected test error (\ref{eq:prerr}) under the aforementioned three settings. To make  fair comparisons, we let $m=n$ in the Random-$X$ setting.
However, the theoretical results that we have developed for  expected test errors under the Random-$X$ setting are valid for the test data
with general size $m$.  In this section, we also assume that the true covariance regression model is $\Cov[Y_\alpha|\mF_{\mathcal{X}}]=C(X_\alpha)$ for both training and test data (i.e., $\alpha\in\{1,\cdots,n\}\cup\{01,\cdots,0n\}$). Note that $C(\cdot)$ does not necessarily have the same parametric structure as $C_\beta(\cdot)$, and this allows us to study the impact of model misspecification on the resulting  expected test errors.
With a slight abuse of notation,
we incorporate the training and test data together and denote $\mF_{\mathcal{X}}=\sigma\langle X_1,\cdots,$ $X_n,X_{01},\cdots,X_{0n}\rangle$.  In what follows, we introduce the  main theoretical results with regard to the expected  test errors. Additional theoretical results and remarks are provided in Section \ref{sec:adS3} of the supplementary material.

\subsection{Bias-variance decomposition and optimism}

\label{sec:b-v}
To study  the expected test error $\Err$ in (\ref{eq:prerr}), we introduce a technical assumption below. This assumption is not only applied for both the training and test data,
but also
allows us to compare the bias-variance decompositions of $\Err$
across the Fixed-$X$, Random-$X$ and Same-$X$ settings.

\setcounter{as}{5}


\begin{as}\label{assump:1} Assume that the data are generated by
$
Y_{\alpha}=C^{1/2}(X_\alpha)\mE_{\alpha}$
$\textrm{ for all }\alpha=1,\cdots,n$ and $01,\cdots,0n$, $\mE_{\alpha}=(\varepsilon_{\alpha 1}, \cdots, \varepsilon_{\alpha p})^\top
$ is independent of $\mF_X$,
 where
 $\varepsilon_{\alpha j}$s are  i.i.d. in both $\alpha$ and $j$ with $\rmE \,(\varepsilon_{\alpha j})=0$, $\rmE\, (\varepsilon_{\alpha j}^2)=1$,  $\rmE\, (\varepsilon_{\alpha j}^4)=\mu^{(4)}$,  $\rmE|\varepsilon_{\alpha j}|^{4+\eta_1}<\infty$ for   $j=1,\cdots,p$, and $\eta_1>0$ is defined in Assumption \ref{as:1}.

\end{as}

A condition similar to Assumption \ref{assump:1} is also given in  \citet{Rosset2020} for establishing the  model assessment theory of  classical mean regression models.  
 Assumption \ref{assump:1} is required for studying
 the bias-variance decomposition under the Fixed-$X$, Same-$X$, 
and Random-$X$ settings.   
To this end,  we first obtain the closed form of the conditional covariance $\bmV(X_\alpha)=\Cov[\bmY_\alpha|X_\alpha]$ for $\alpha\in\{1,\cdots,n\}\cup\{01,\cdots,0n\}$; see Proposition \ref{pn:ccov} of Section \ref{sec:adS3} in the supplementary material. 
In order to tackle the technical challenges and 
complexities of establishing the model assessment theory due to the complicated dependence structure $\Cov[\bmY_\alpha|X_\alpha]$, the training error (\ref{eq:traininge}) and test error (\ref{eq:teste})  defined earlier have  incorporated the  general weight matrix $\bmW(X_\alpha)$. This is a new idea compared to the training error and test error considered by \citet{Rosset2020}, who used the identity weight matrix.
Then, let  $\nu^2(X_\alpha)=\tr\{\bmW(X_\alpha)\bmV(X_\alpha)\}$,
which is non-negative. In the bias-variance decomposition, we let  $\mB$ and $\mV$ correspondingly represent  the bias and variance measures of  the estimator $\hat{\widetilde \bmC}(X_\alpha)$ of $\widetilde \bmC(X_\alpha)$, where $\widetilde \bmC(X_\alpha)=\bmW^{1/2}(X_\alpha) \bmC(X_\alpha)$ with $\bmC (X_\alpha)=\rmvec\, C(X_\alpha)$,  and $\hat{\widetilde \bmC}(X_\alpha)=\bmW^{1/2}(X_\alpha)\hat \bmC(X_\alpha)$. We further denote  $\widetilde \bmY_\alpha=\bmW^{1/2}(X_\alpha)\bmY_{\alpha}$. Finally, we evaluate the expected test errors $\Err$ in (\ref{eq:prerr})  under the Fixed-$X$, Random-$X$ and Same-$X$ settings, and obtain $\ErrF$, $\ErrR$ and $\ErrS$, respectively; see Proposition \ref{pn:Err} in Section \ref{sec:adS3} of the supplementary material.  Based on the results in  Propositions  \ref{pn:ccov} and \ref{pn:Err}, we next show the bias-variance decompositions for $\ErrF$, $\ErrR$ and $\ErrS$, respectively.


\begin{tm}\label{tm:P1} Suppose Assumption \ref{assump:1} is satisfied.
(i) Under the Fixed-$X$ setting, we have  $\ErrF={n^{-1}}\sum_{i=1}^n$ $ \nu^2(x_i)+\mB_{\ErrF}+\mV_{\ErrF}$, where
\[
\mB_{\ErrF}=\frac{1}{n}\sum_{i=1}^n\laak \rmE\,\hat \wbmC(x_i)- \wbmC(x_i)\raak_2^2\textrm{ and }\mV_{\ErrF}=\frac{1}{n}\sum_{i=1}^n\tr\lbk  \Cov\lsk \hat \wbmC(x_i) \rsk \rbk;
\]
(ii) Under the Random-$X$ setting, we have $\ErrR= \rmE\,  \big\{\nu^2(X_{1})\big\}+\mB_{\ErrR}+\mV_{\ErrR}$,  where
\[
\mB_{\ErrR}=\rmE\laak \rmE\lmk \hat \wbmC(X_{01})\Big{|}\mbX,X_{01}\rmk- \wbmC(X_{01})\raak_2^2\textrm{ and }\mV_{\ErrR}=\rmE\,\tr\lbk  \Cov\lmk \hat \wbmC(X_{01})\Big{|}\mbX,X_{01}\rmk \rbk;
\]
(iii) Under the Same-$X$ setting, we have $\ErrS= \rmE \, \big\{\nu^2(X_{1})\big\}+\mB_{\ErrS}+\mV_{\ErrS}$, where
\[
\mB_{\ErrS}=\rmE \laak \rmE\lmk \hat \wbmC(X_{1})\Big{|}\mbX\rmk- \wbmC(X_{1})\raak_2^2\textrm{ and }\mV_{\ErrS}=\rmE\,\tr\lbk  \Cov\lmk \hat \wbmC(X_{1}) \Big{|}\mbX\rmk \rbk.
\]
\end{tm}

\noindent
It is worth noting that $\mB_{\ErrF}$ and  $\mV_{\ErrF}$ in (i) represent the  mean squared bias and the average trace of $\textrm{Cov}\big(\hat{\widetilde \bmC}(x_i)\big)$, respectively.
This  interpretation is also applicable to ($\mB_{\ErrR}$, $\mV_{\ErrR}$) in (ii)
and  ($\mB_{\ErrS}$, $\mV_{\ErrS}$) in (iii). For the sake of simplicity, we refer to $\mB$ and $\mV$ as ``bias'' and ``variance'', respectively.

To assess the difference between the Random-$X$ and Same-$X$ settings, we employ Theorem
 \ref{tm:P1} and obtain that
\be\label{eq:R-S}
\ErrR-\ErrS=\mB^{(d)}+\mV^{(d)},
\ee
where $\mB^{(d)}=\mB_{\ErrR}-\mB_{\ErrS}$ and $\mV^{(d)}=\mV_{\ErrR}-\mV_{\ErrS}$
represent the bias and variance differences, respectively.
If $\mB^{(d)}+\mV^{(d)}\geq 0$,
then  the Random-$X$ error,  $\ErrR$, is no smaller than
the Same-$X$  error, $\ErrS$.
The  difference between $\ErrR$ and $\ErrS$ in (\ref{eq:R-S}) allows us to
make comparisons between  $\ErrR$ and $\ErrF$, which is our main interest and is discussed below.

We first adapt the concept of optimism  from classical mean regression (see, e.g., \citealp{Hastie2009}), and then define the optimism in covariance regression as follows: \be\label{eq:Opt}
\Opt=\rmE\lbk  \Te(\mbX_0,\mbY_0,\mbX,\mbY)-\Tr(\mbX,\mbY) \rbk=\Err- \rmE\lbk \Tr(\mbX,\mbY)  \rbk,
\ee
which is the mean of the difference between the test error and training error.
We next evaluate $\Opt$ under the Fixed-$X$, Same-$X$ and  Random-$X$  settings and obtain $\OptF$, $\OptS$ and $\OptR$, respectively.


\begin{pn}\label{pn:Opt} Suppose Assumption \ref{assump:1} is satisfied. (i) Under the Fixed-$X$ setting,  we have $\OptF={2}n^{-1}\sum_{i=1}^n\tr\,\Cov\big( \widetilde \bmY_i,\hat\wbmC(x_i)\big)$; (ii) under the Same-$X$ setting, we obtain $\OptS$ $={2}n^{-1}\sum_{i=1}^n\tr\,\rmE\,\Cov\big[ \widetilde \bmY_i,\hat\wbmC(X_i)\big{|}\mbX\big]$; and
(iii) under the Random-$X$ setting, we have
$
\OptR=\OptS+\mB^{(d)}+\mV^{(d)}
$.
\end{pn}

Using  Proposition \ref{pn:Opt}, we can further explore the relationships between the estimators of
 $\ErrF$, $\ErrS$ and $\ErrR$.
Let $U(\mbX,\mbY)$ be an unbiased estimator of $\OptF$ under the Fixed-$X$ setting, and a typical
example of the unbiased estimator of $\OptF$ is presented in Section \ref{sec:non-negative}.  Based on (\ref{eq:Opt}), $\ErrF$ can be naturally estimated by $\widehat\ErrF=\Tr(\mbX,\mbY)+U(\mbX,\mbY)$, which is also the unbiased estimator of $\ErrS$ demonstrated below.


\begin{pn} \label{pn:OptS} Suppose Assumption \ref{assump:1} is satisfied. (i) $U(\mbX,\mbY)$ is an unbiased estimator of $\OptS$.
In addition, the unbiased estimator of  the expected  test error in (\ref{eq:prerr}) is the same under the Fixed-$X$ and Same-$X$ settings. Specifically,
\be\label{eq:hatErrF}
\widehat\ErrF=\widehat\ErrS=\Tr(\mbX,\mbY) +U(\mbX,\mbY).
\ee
 (ii) If $\OptF$  does not depend on $x_1,\cdots,x_n$ under the Fixed-$X$ setting,
then it is equal to $\OptS$ under the Same-$X$ setting.
\end{pn}

Let $\widehat\mB^{(d)}$ and $\widehat\mV^{(d)}$ be the estimators of $\mB^{(d)}$ and  $\mV^{(d)}$, respectively.  An example of the estimators of  $\mB^{(d)}$ and  $\mV^{(d)}$ is discussed in Section \ref{sec:non-negative}. By  Propositions \ref{pn:Opt} -- \ref{pn:OptS} and  (\ref{eq:Opt}),
we finally obtain that
\[
\widehat\OptR=U(\mbX,\mbY)+\widehat\mB^{(d)}+\widehat\mV^{(d)}\textrm{, and }
\]
\be\label{eq:hatErrR}
\widehat\ErrR=\Tr(\mbX,\mbY)+\widehat\OptR=\Tr(\mbX,\mbY) +U(\mbX,\mbY)+\widehat\mB^{(d)}+\widehat\mV^{(d)},
\ee
\noindent
In sum, equations (\ref{eq:hatErrF}) and (\ref{eq:hatErrR}) clearly indicate the
components that are due to, first, the difference $\mB^{(d)}$ between the biases in the estimators of $\ErrF$ and $\ErrR$ and, second, the difference $\mV^{(d)}$ between the variances in  the estimators of $\ErrF$ and $\ErrR$.

In classical mean regression models,
$\widehat\ErrF$  can refer to  Mallows' Cp \citep{Mallows1973}, which is commonly used as a model  assessment criterion.  Analogously, $\widehat\ErrR$ can be used as a model assessment criterion for random covariance regressions, and that will be discussed in the next subsection.

\subsection{Mallows' Cp in linear covariance regression models}

\label{sec:non-negative}

In this section, we mainly focus on model assessment based on fitting linear covariance regression model (\ref{eq:lcov}) with the estimator
$\hat\beta_{\WLS}$ defined in (\ref{eq:GLS}). However, the true covariance regression model $\Cov[Y_\alpha|\mF_{\mathcal{X}}]=C(X_\alpha)$ for both training and test data (i.e., $\alpha\in\{1,\cdots,n\}\cup\{01,\cdots,0n\}$) is not necessarily correctly specified by
(\ref{eq:lcov}).

Note that in model (\ref{eq:lcov}), the $X_i^{(k)}$s are $p\times p$  matrices for $i=1,\cdots,n$ and $k=1,\cdots,K$.  Denote $\bmX_i=\big(\rmvec(X_i^{(1)}),\cdots, \rmvec(X_i^{(K)})\big)\in\mathbb{R}^{p^2\times K}$.
In addition, let $X_i=(X_i^{(1)},\cdots,X_i^{(K)})$, which was defined in  the first example  in Section \ref{sec:review}.
Under the Fixed-$X$ setting, consider $X_i=x_i=(x_i^{(1)},\cdots,x_i^{(K)})$ and
$\bmx_i=\big(\rmvec(x_i^{(1)}),\cdots,\rmvec(x_i^{(K)})\big) \in\mathbb{R}^{p^2\times K}$.   
 Let $x=\{x_1,\cdots,x_n\}$ and $\mbW(x)=\textrm{diag}\{\bmW(x_1),\cdots,\bmW(x_n)\}$ be a $np^2\times np^2$ matrix whose diagonal blocks are $\bmW(x_1),\cdots,\bmW(x_n)$ and the other elements are zeros.  Similarly, let $\mbV(x)=\textrm{diag}\{ \bmV(x_1) ,\cdots, \bmV(x_n) \}\in \mathbb{R}^{np^2\times np^2}$,  and the detailed expression of  $\bmV(\cdot) \in\mathbb{R}^{p^2\times p^2}$ is  presented in Proposition \ref{pn:ccov} of Section \ref{sec:adS3} of the supplementary material. Finally, let $\widetilde \bmX_i=\bmW^{1/2}(X_i)\bmX_i \in\mathbb{R}^{p^2\times K}$, $\widetilde \bmx_i=\bmW^{1/2}(x_i)\bmx_i\in\mathbb{R}^{p^2\times K}$, $
\widetilde\bmX=(\widetilde\bmX_1^\top,\cdots,\widetilde\bmX_n^\top)^\top \in\mathbb{R}^{np^2\times K}$, $
\widetilde\bmx=(\widetilde\bmx_1^\top,\cdots,\widetilde\bmx_n^\top)^\top \in\mathbb{R}^{np^2\times K}$,  and $\widetilde\mbV(x)=\mbW^{1/2}(x)\mbV(x)\mbW^{1/2}(x)\in \mathbb{R}^{np^2\times np^2}$.  Analogously, for $\alpha\in\{01,\cdots,0n\}$, we denote $\bmX_\alpha=\big(\rmvec(X_\alpha^{(1)}),\cdots, \rmvec(X_\alpha^{(K)})\big)\in\mathbb{R}^{p^2\times K}$, $X_\alpha=(X_\alpha^{(1)},\cdots,$ $X_\alpha^{(K)})$, and $\widetilde \bmX_\alpha=\bmW^{1/2}(X_\alpha)\bmX_\alpha\in\mathbb{R}^{p^2\times K}$. The above notation will be used in the following  proposition and theorem.


\begin{pn}\label{pn:OptOLS} Suppose Assumption \ref{assump:1} is satisfied. In addition,  consider fitting a linear covariance regression model
(\ref{eq:lcov}) with  $\hat\beta=\hat\beta_{\WLS}$ defined in (\ref{eq:GLS}).
Under the Fixed-$X$ setting,  we have that $
\OptF={2}{n^{-1}}\tr\big\{ $ $  ( \widetilde\bmx^\top \widetilde\bmx
)^{-1}\widetilde\bmx^\top \widetilde\mbV(x) \widetilde \bmx\big\}.
$
In addition, under the Same-$X$ setting,  we have that $
\OptS={2}{n^{-1}}\rmE\,\tr\big\{   ( \widetilde\bmX^\top \widetilde\bmX
)^{-1}\widetilde\bmX^\top \widetilde\mbV(\mbX) \widetilde \bmX\big\}.
$
\end{pn}

As discussed in  Remark \ref{re:m1}, set $\bmW(X_\alpha)=\bmV^{-}(X_\alpha)$. Then
 Proposition \ref{pn:OptOLS} leads to
\be\label{eq:OptFS}
\OptF=\OptS=\frac{2K}{n},
\ee
which is a special case of the general result, $\OptF=\OptS$, as demonstrated in Proposition \ref{pn:OptS}.
Note that  ${2K}/{n}$ is a constant and does not need to be estimated.
 In addition, consider the case that  $\bmV(\cdot)$ is known. This, together with Proposition \ref{pn:OptOLS}, yields an unbiased estimator for both $\OptF$ and $\OptS$
as given below.
\be\label{eq:ux}
U(\mbX,\mbY)=\frac{2}{n}\tr\big\{   ( \widetilde\bmX^\top \widetilde\bmX
)^{-1}\widetilde\bmX^\top \widetilde\mbV(\mbX) \widetilde \bmX\big\}\triangleq U(\mbX).
\ee
Note that this estimator does not depend on $\mbY$.
Using the above result and Proposition \ref{pn:OptS}, we define Mallows' Cp for linear covariance regression models under the Fixed-$X$ and Same-$X$ settings as follows:
\be\label{eq:Cp}
\Cp\triangleq \widehat\ErrF=\widehat\ErrS=\Tr(\mbX,\mbY) +U(\mbX),
\ee
 which is an unbiased estimator of  $\ErrF$ or $\ErrS$ based on Proposition \ref{pn:OptS}.


\begin{re}\label{re:4}
It is worth noting that the estimator $U(\mbX)$ in (\ref{eq:ux})  depends on
$\bmV(\cdot)$. In  practice, $\bmV(\cdot)$ may be unknown and needs to be estimated. If $C(\cdot)$ is correctly specified by the linear covariance regression model (\ref{eq:lcov}),  Section \ref{subsec:SC} of the supplementary material provides an estimator
$\hat\bmV(\cdot)\triangleq\hat\bmV_{\mbX,\mbY}(\cdot)$ of $\bmV(\cdot)$ based on the training data $(\mbX,\mbY)$. If  we replace $\bmV(\cdot)$ by $\hat\bmV_{\mbX,\mbY}(\cdot)$ in the expression of $U(\mbX)$,  we then obtain the  estimator $\widehat U_{\mbX,\mbY}(\mbX)$. This replacement yields $\Cp$ for the unknown $\bmV(\cdot)$, which is
\be\label{eq:CpE}
\widehat \Cp=\Tr(\mbX,\mbY) +\widehat U_{\mbX,\mbY}(\mbX).
\ee
\end{re}

To  evaluate the estimator of $\ErrR$ and compare it with  $\Cp$ and $\widehat \Cp$, we next employ Proposition \ref{pn:OptOLS} and Theorem \ref{tm:P1} to demonstrate that the
difference between the biases, which is $\mB^{(d)}$, and the difference between the variances, which is $\mV^{(d)}$,   are both nonnegative.


\begin{tm} \label{eq:excess} Suppose Assumption \ref{assump:1} is satisfied. Consider fitting a linear covariance regression model
(\ref{eq:lcov}) with  $\hat\beta=\hat\beta_{\WLS}$ defined in (\ref{eq:GLS}).
We then have  $\mV_{\ErrS}={n^{-1}}\rmE\,\tr\big\{   ( \widetilde\bmX^\top \widetilde\bmX
)^{-1}\widetilde\bmX^\top \widetilde\mbV(\mbX) $ $\widetilde \bmX\big\}=\rmE\, U(\mbX)/2=\OptS/2$,
\bea
\mV_{\ErrR}&=&\frac{1}{n}\tr\lmk \rmE\lbk   \lsk  \widetilde\bmX^\top \widetilde\bmX
\rsk^{-1}\widetilde\bmX^\top \widetilde\mbV(\mbX) \widetilde \bmX  \lsk  \widetilde\bmX^\top \widetilde\bmX
\rsk^{-1}\rbk\rmE \lsk \widetilde\bmX^\top \widetilde\bmX \rsk\rmk\nn\\
&=& \rmE\,\tr\lbk\lsk \widetilde\bmX^\top \widetilde\bmX
\rsk^{-1}\widetilde\bmX^\top\widetilde\mbV(\mbX) \widetilde \bmX \lsk \widetilde\bmX^\top \widetilde\bmX
\rsk^{-1} \widetilde \bmX_{01}^\top \widetilde\bmX_{01}\rbk
,\textrm{ and }\label{eq:verrro}
\eea
$\mB^{(d)}=\mB_{\ErrR}-\mB_{\ErrS}\geq 0$. Furthermore, set $\bmW(X_\alpha)=\bmV^{-}(X_\alpha)$. We   obtain $\mV^{(d)}=\mV_{\ErrR}-\mV_{\ErrS}\geq 0$, and hence  $\ErrR\geq \ErrS$.
\end{tm}


\begin{re}\label{re:5}
Suppose that $C(\cdot)$ is correctly specified by the linear covariance regression model (\ref{eq:lcov}).
Then we can verify that $\mB_{\ErrR}=\mB_{\ErrS}=0$, which implies that the bias difference is $\mB^{(d)}=0$.
In addition,  consider the ideal case that $\hat\beta=\beta_0$. We then have
 $\hat{\widetilde \bmC}(X_\alpha)=\widetilde \bmC(X_\alpha)$. This, in conjunction with   Proposition \ref{pn:Err} in Section \ref{sec:adS3} of the supplementary material, leads to $\ErrR=\ErrS$.
Moreover, Theorem \ref{tm:3} shows that $\hat\beta_{\WLS}$ is the consistent
estimator of $\beta_0$ as $\max\{n,p\}\to\infty$. Even though $\hat\beta_{\WLS}\neq\beta_0$, $\hat\beta_{\WLS}$ converges to $\beta_0$ with probability tending to 1. Using this result and  Proposition \ref{pn:Err}, one can expect that, as $\max\{n,p\}\to\infty$, the difference between $\ErrR$ and $\ErrS$ tends to be small under the correctly specified model.

\end{re}


 Note that, under $\bmW(X_\alpha)=\bmV^{-}(X_\alpha)$,
(\ref{eq:Opt})  and (\ref{eq:OptFS}) imply that $\ErrS=\rmE\{\Tr(\mbX,\mbY)\}+2K/n$ and $\ErrF=\rmE\{\Tr(x,\mbY)\}+2K/n$.
Hence, $\ErrS$ and $\ErrF$ are different.
In practice, however, $\ErrS$ and $\ErrF$ have the same estimator $\Cp$ in (\ref{eq:Cp}), and it is different from the $\ErrR$ estimator, $\widehat\ErrR$,  in (\ref{eq:hatErrR}). From Theorem \ref{eq:excess}, one can expect that the estimators of the difference in biases and the difference in variances satisfy $\widehat \mB^{(d)}\geq 0$ and  $\widehat\mV^{(d)}\geq 0$, respectively, which leads to $\widehat\ErrR\geq \Cp$. We next propose  two estimation
approaches to obtain $\widehat\ErrR$ explicitly.

\subsubsection{Estimation of $\ErrR$ with $C(\cdot)$  correctly specified in  model (\ref{eq:lcov})}

\label{sec:RCp}

Assume that $C(\cdot)$ is correctly specified by the linear covariance regression model (\ref{eq:lcov}).
By  Remark \ref{re:5}, we have $\mB^{(d)}=\mB_{\ErrR}=\mB_{\ErrS}=0$.
Naturally, their corresponding estimators are  $\widehat\mB^{(d)}=\widehat\mB_{\ErrR}=\widehat\mB_{\ErrS}=0$.
Hence, based on (\ref{eq:R-S}) and (\ref{eq:Cp}),
it suffices to estimate $\mV^{(d)}=\mV_{\ErrR}-\mV_{\ErrS}$ and then obtain the estimator of $\ErrR$.
We next discuss the  estimation of $\mV_{\ErrS}$ and  $\mV_{\ErrR}$, separately.

Under the assumption that $\bmV(\cdot)$ is known,
$U(\mbX)$ in (\ref{eq:ux}) is an unbiased estimator of $\OptS$.
 In addition by Theorem \ref{eq:excess},  we obtain that
$\widetilde \mV_{\ErrS}=U(\mbX)/2$ is an unbiased estimator of $\mV_{\ErrS}$.
 As for the estimation of $\mV_{\ErrR}$ in (\ref{eq:verrro})  of Theorem \ref{eq:excess}, we only have the training data $\mbX$ and  $\widetilde\bmX$, and the test data $\widetilde\bmX_{01}$ is not available.
Hence, we consider using 
\be\label{eq:almost}
\rmE\,\tr\lbk\lsk \widetilde\bmX^{(-i)}{}^\top \widetilde\bmX^{(-i)}
\rsk^{-1}\widetilde\bmX^{(-i)}{}^\top\widetilde\mbV(\mbX^{(-i)}) \widetilde \bmX^{(-i)} \lsk \widetilde\bmX^{(-i)}{}^\top \widetilde\bmX^{(-i)}
\rsk^{-1} \widetilde \bmX_{i}^\top \widetilde\bmX_{i}\rbk,
\ee
to approximate (\ref{eq:verrro}), where $
\widetilde\bmX^{(-i)}=(\widetilde\bmX_1^\top,\cdots,\widetilde\bmX_{i-1}^\top,\widetilde\bmX_{i+1}^\top,\cdots,\widetilde\bmX_n^\top)^\top$ and  $\mbX^{(-i)}=\{X_1,\cdots,X_{i-1},X_{i+1},\cdots, X_n\}$ for $i=1,\cdots,n$. Compared to the training data  $\mbX$ of size $n$ in (\ref{eq:verrro}),  $\mbX^{(-i)}$ in (\ref{eq:almost}) is of size $n-1$. As a consequence, one  expects that the approximation (\ref{eq:almost}) is close to $\mV_{\ErrR}$ as $n$ gets large. Thus, it is natural to consider a cross-validation estimator of  $\mV_{\ErrR}$,
\be\label{pn:CV0}
\widetilde \mV_{\ErrR}=\frac{1}{n} \sum_{i=1}^n\tr\lbk\lsk \widetilde\bmX^{(-i)}{}^\top \widetilde\bmX^{(-i)}
\rsk^{-1}\widetilde\bmX^{(-i)}{}^\top\widetilde\mbV(\mbX^{(-i)}) \widetilde \bmX^{(-i)} \lsk \widetilde\bmX^{(-i)}{}^\top \widetilde\bmX^{(-i)}
\rsk^{-1} \widetilde \bmX_{i}^\top \widetilde\bmX_{i}\rbk,
\ee
which is an unbiased estimator of (\ref{eq:almost}). Based on the estimators $\widetilde\mV_{\ErrR}$ and $\widetilde\mV_{\ErrS}$,
 we then obtain an estimator of $\mV^{(d)}$, $\widetilde\mV^{(d)}=\widetilde\mV_{\ErrR}-\widetilde\mV_{\ErrS}$.

The equation $\widetilde\mV^{(d)}=\widetilde\mV_{\ErrR}-\widetilde\mV_{\ErrS}$,  along with (\ref{eq:R-S}), (\ref{eq:Cp}) and the fact that
$\widetilde \mV_{\ErrS}=U(\mbX)/2$, yields an estimator of $\ErrR$, which is
\be \label{pn:CVR}
\RCp\triangleq\Cp+\widetilde\mV^{(d)}=\Tr(\mbX,\mbY) +\frac{1}{2}U(\mbX)+\widetilde\mV_{\ErrR}.
\ee
This estimator is essentially  Mallows' Cp adjusted by $\widetilde\mV^{(d)}$,
the estimator of the difference in variances, which is induced by the Random-$X$ setting. Thus, we name it Random Cp ($\RCp$) for covariance regression.
Note that $\bmV(\cdot)$ may be unknown in practice.
 Analogously to  Remark \ref{re:4},
we can replace $\bmV(\cdot)$ in (\ref{pn:CV0}) by its estimator $\hat \bmV(\cdot)$ defined in Section \ref{subsec:SC} of the supplementary material. Then, we obtain the estimator of $\mV_{\ErrR}$ for the unknown $\bmV(\cdot)$, and denote it by
$\widehat\mV_{\ErrR}$. This, together with the estimator of $U(\mbX)$ in Remark \ref{re:4},
namely  $\widehat U_{\mbX,\mbY}(\mbX)$,
 yields the estimator of $\ErrR$:
\be \label{pn:CVRE}
\widehat \RCp=\Tr(\mbX,\mbY) +\frac{1}{2}\widehat U_{\mbX,\mbY}(\mbX)+\widehat\mV_{\ErrR}.
\ee
Based on Theorem \ref{eq:excess}, we expect that $\widehat \RCp$ is larger than $\widehat \Cp$.
It is worth noting that either $\widehat \RCp$ or $\widehat \Cp$ can be used to compare
a set of covariance regression models in terms of  their expected test errors.
 Under the Random-$X$ setting, however,  $\widehat \RCp$ should be a better choice since $\widehat \Cp$ can underestimate expected  test errors. Note that  $\RCp$ and $\widehat \RCp$ use the results from  Proposition \ref{pn:OptOLS} and Theorem \ref{eq:excess}
obtained by  fitting a linear covariance regression model (\ref{eq:lcov}). Hence, they are only applicable for the linear model.
The following subsection discusses the second approach to estimating $\ErrR$.

\subsubsection{Estimation of $\ErrR$ with $C(\cdot)$ misspecified in  model (\ref{eq:lcov})}

\label{sec:OCV}

Recall that before equation (\ref{eq:traininge}), we denoted $\hat \bmC(\cdot)=\bmC_{\hat\beta}(\cdot)$.
 If $C(\cdot)$ is misspecified by the linear covariance regression model (\ref{eq:lcov}),   $\hat \bmC(\cdot)=\bmC_{\hat\beta}(\cdot)$ is a biased estimator of $\bmC(\cdot)$. Hence,
 $\mB_{\ErrR}\neq 0$ and  $\mB_{\ErrS}\neq 0$, and $\widehat \RCp$ is not an appropriate estimator of $\ErrR$.
Since the true covariance function $C(\cdot)$ is unknown, we cannot estimate $\mB_{\ErrR}$ and $\mB_{\ErrS}$  directly based on Theorem \ref{tm:P1}.  Instead, we use the result $
\ErrR= \rmE\|  \widetilde \bmY_{01}-\hat \wbmC(X_{01})\|_2^2
$ in Proposition \ref{pn:Err}  of Section \ref{sec:adS3} of the supplementary material,
  adopt a similar approach to obtaining the estimator (\ref{pn:CV0}),
 and then find a cross-validation estimator. Specifically,  let ${\hat \wbmC}^{(-i)}(\cdot)$ be the WLS estimator  obtained from all training data $(\mbX,\mbY)$ except for $(X_i,Y_i)$.
Then, calculate the ordinary cross-validation ($\OCV$) error,
\be\label{eq:nocv}
\OCV=\frac{1}{n}\sum_{i=1}^n \laak  \widetilde \bmY_{i}-\hat \wbmC^{(-i)}(X_{i}) \raak_2^2,
\ee
which is an approximately unbiased estimator of $\ErrR$.

The advantage of using $\OCV$ is that it can estimate $\ErrR$
for a general and possibly nonlinear $C_\beta(\cdot)$. In addition, it does not require the true covariance function $C(\cdot)$ to be correctly specified. If $C(\cdot)$ is correctly specified by the linear covariance regression model (\ref{eq:lcov}), then either $\OCV$ or $\widehat\RCp$ can be used to estimate $\ErrR$. Under the random covariance regression models,  we compare the performance of
the statistics $\Cp$, $\widehat\Cp$, $\RCp$, $\widehat\RCp$ and $\OCV$  via  simulation studies in Section \ref{sec:40} below.



\section{Simulation Studies}
\label{sec:40}



In this  section, we conduct simulation studies to investigate the finite sample performance of parameter estimations and  expected test errors in covariance regression models.   All the tables and  figures in this section, except Tables \ref{tb:EC1m-1} -- \ref{tb:EC1m-2},  are provided in Section \ref{sec:tbfig} of the supplementary material to save space.

In these studies, for both training data $\alpha\in\{1,\cdots,n\}$ and test data $\alpha\in\{01,\cdots,0n\}$, the response vectors $Y_\alpha$ are simulated by
$Y_{\alpha}=C^{1/2}_{\beta_0}(X_\alpha)\mE_{\alpha}$, where the components of the random error $\mE_{\alpha}$ are independent and identically generated from the standard normal distribution and the mixture of normal distributions $0.9N(0,5/9)+0.1N(0,5)$, respectively. In addition, the following two different models are considered
for $C_{\beta_0}(X_\alpha)$ with $X_\alpha=(X_\alpha^{(2)},\cdots,X_\alpha^{(K)})$
and $K=5$.

\textbf{Model A.} Let $C_{\beta_0}(X_\alpha)=\beta_0^{(1)}I_p+\sum_{k=2}^{5}\beta_0^{(k)}X_\alpha^{(k)}$,
where the $X_\alpha^{(k)}$s are $p\times p$  matrices, the diagonal elements of $X_\alpha^{(k)}$ are zeros for $k=2,\cdots,5$, the off-diagonal elements of $X_\alpha^{(2)}$ are independent and identically generated from the Bernoulli distribution with probability $p^{-1}$, and the off-diagonal elements of $X_\alpha^{(k)}=(X_{\alpha,j_1j_2}^{(k)})_{p\times p}$ for $k=3,4,5$ are $X_{\alpha,j_1j_2}^{(k)}=\exp\{-(d_{\alpha,j_1j_2}^{(k)})^2\}$ with $d_{\alpha,j_1j_2}^{(k)}$  independent and identically generated from the uniform distribution $U(0,p^{1+(k-3)/6})$.

\textbf{Model B.} Let
 $C_{\beta_0}(X_\alpha)=\beta_0^{(1)}I_p+\sum_{k=2}^5 \beta_0^{(k)} X_\alpha^{(k)}X_\alpha^{(k)\top}$, where the $X_\alpha^{(k)}$s are $p\times 1$  vectors, and the components of $X_\alpha^{(k)}$  are independent and identically generated from the normal distribution $N(0,p^{-1})$, for $k=2,\cdots, 5$.

Note that Models A and B, respectively, correspond to the first example and a variant of the  third example of  covariance regression in  Section \ref{sec:review}. In both models, we set $\beta_0^{(1)}=5$ and $\beta_0^{(k)}=1.4$ for $k=2,\cdots,5$.
In addition, we consider three sample sizes $n=25$, $50$, and $100$, and two dimensions of the response vector $p=5$ and $25$.
For each of the above models, a total of 1,000 simulation experiments
are conducted.
 For the Fixed-$X$ setting, we set $X_\alpha$ with $\alpha=1,\cdots, n$ and $01,\cdots,0n$ such that $X_{0i}=X_i$, for $i=1,\cdots,n$, and we only use one realization of $\{X_1=x_1,\cdots,X_n=x_n\}$ in all 1,000 simulation experiments. In contrast,
for the Random-$X$ setting,  the $X_\alpha$s are independent and identically generated
for $\alpha=1,\cdots, n$ and $01,\cdots,0n$ across 1,000 realizations.
Accordingly, these processes for generating $X_\alpha$ allow us to make a fair comparison between the Fixed-$X$ and Random-$X$ settings.

We first assess the performance of the constrained QMLE ($\hat\beta_{\QMLE}^c$) and  constrained WLS estimators based on the training data $i=1,\cdots,n$, where the constrained estimators are shown to be asymptotically equivalent to the unconstrained estimators in Remark \ref{re:positived}. Among the constrained WLS estimators, we consider the constrained OLS estimator $\hat\beta_{\OLS}^c$  by letting $\bmW(X_i)=I_{p^2}$, the constrained GLS estimator $\hat\beta_{\GLS}^c$ by letting $\bmW(X_i)= \bmV^-(X_i)$,  and the constrained FGLS estimator $\hat\beta_{\FGLS}^c$ by letting $\bmW(X_i)=\hat \bmV^-(X_i)$, where the last two estimators  can be obtained via Remark \ref{re:m1} and Remark \ref{re:positived}.  We also compare the constrained OLS estimate to the unconstrained OLS estimate. The numerical performance of the unconstrained estimate and some related discussion are given in
Section \ref{sec:Ssimu} of the supplementary material.

\begin{table}[htbp!]

\caption{The average bias (BIAS),  standard deviation (SD), and  root mean
squared error (RMSE) of the  constrained QMLE and OLS estimates obtained from Model A with the mixture of normal random errors under the Fixed-$X$ and Random-$X$ settings.}
%
\label{tb:EC1m-1}
\begin{center}
\scalebox{0.99}{
\begin{tabular}{c|c|rrrrr|rrrrr}
  \hline
  & & \multicolumn{5}{c|}{$p=5$}     &    \multicolumn{5}{c}{$p=25$} \\
\cline{3-7}\cline{8-12}
 $n$  &\multicolumn{1}{c|}{}& \multicolumn{1}{c}{$\hat \beta_{1}$} & \multicolumn{1}{c}{$\hat \beta_{2}$} & \multicolumn{1}{c}{$\hat \beta_{3}$} & \multicolumn{1}{c}{$\hat \beta_{4}$} & \multicolumn{1}{c|}{$\hat \beta_{5}$} & \multicolumn{1}{c}{$\hat \beta_{1}$} & \multicolumn{1}{c}{$\hat \beta_{2}$} & \multicolumn{1}{c}{$\hat \beta_{3}$} & \multicolumn{1}{c}{$\hat \beta_{4}$} & \multicolumn{1}{c}{$\hat \beta_{5}$}   \\
   \hline

   \multicolumn{12}{c}{QMLE \& Fixed-$X$ Setting}\\
   \hline
25 & BIAS & 0.34 & 0.01 & -0.23 & -0.19 & -0.15 & 0.28 & 0.02 & 0.04 & -0.04 & -0.06 \\
 & SD & 0.96 & 0.72 & 0.81 & 0.84 & 0.90 & 0.41 & 0.20 & 0.28 & 0.32 & 0.47 \\
 & RMSE & 1.02 & 0.72 & 0.84 & 0.86 & 0.92 & 0.50 & 0.20 & 0.28 & 0.33 & 0.47 \\
   \hline
  50 & BIAS & 0.36 & -0.05 & 0.06 & -0.04 & -0.18 & 0.23 & 0.00 & 0.01 & 0.05 & 0.05 \\
 & SD & 0.67 & 0.39 & 0.56 & 0.57 & 0.57 & 0.30 & 0.16 & 0.20 & 0.26 & 0.32 \\
 & RMSE & 0.76 & 0.39 & 0.56 & 0.57 & 0.60 & 0.38 & 0.16 & 0.20 & 0.27 & 0.32 \\
   \hline
  100 & BIAS & 0.32 & 0.08 & 0.08 & -0.10 & -0.05 & 0.23 & 0.04 & 0.04 & -0.00 & -0.00 \\
 & SD & 0.47 & 0.32 & 0.37 & 0.37 & 0.41 & 0.19 & 0.11 & 0.13 & 0.15 & 0.19 \\
 & RMSE & 0.57 & 0.33 & 0.38 & 0.38 & 0.41 & 0.30 & 0.12 & 0.13 & 0.15 & 0.19 \\
\hline
\multicolumn{12}{c}{QMLE \& Random-$X$ Setting}\\
   \hline
25 & BIAS & 0.11 & -0.12 & -0.11 & -0.13 & -0.11 & 0.15 & -0.02 & -0.00 & -0.01 & -0.01 \\
 & SD & 1.02 & 0.58 & 0.70 & 0.80 & 0.91 & 0.47 & 0.21 & 0.26 & 0.34 & 0.42 \\
 & RMSE & 1.03 & 0.60 & 0.71 & 0.81 & 0.92 & 0.49 & 0.22 & 0.26 & 0.34 & 0.42 \\
   \hline
  50 & BIAS & 0.17 & -0.05 & -0.03 & -0.06 & -0.07 & 0.18 & 0.01 & 0.01 & 0.03 & 0.04 \\
 & SD & 0.72 & 0.39 & 0.53 & 0.57 & 0.60 & 0.33 & 0.16 & 0.19 & 0.26 & 0.31 \\
 & RMSE & 0.74 & 0.39 & 0.53 & 0.57 & 0.60 & 0.38 & 0.16 & 0.19 & 0.26 & 0.31 \\
   \hline
  100 & BIAS & 0.17 & -0.02 & 0.01 & -0.02 & -0.03 & 0.20 & 0.02 & 0.02 & 0.02 & 0.03 \\
 & SD & 0.48 & 0.28 & 0.36 & 0.38 & 0.40 & 0.20 & 0.11 & 0.13 & 0.17 & 0.23 \\
& RMSE & 0.51 & 0.28 & 0.36 & 0.38 & 0.40 & 0.28 & 0.11 & 0.13 & 0.17 & 0.23 \\
\hline

\multicolumn{12}{c}{OLS \& Fixed-$X$ Setting}\\
   \hline
25 & BIAS & -0.11 & -0.33 & -0.54 & -0.43 & -0.45 & -0.10 & -0.18 & -0.14 & -0.26 & -0.24 \\
 & SD & 1.27 & 0.76 & 0.85 & 0.86 & 0.96 & 0.51 & 0.28 & 0.35 & 0.43 & 0.59 \\
 & RMSE & 1.27 & 0.83 & 1.01 & 0.96 & 1.06 & 0.52 & 0.33 & 0.38 & 0.50 & 0.64 \\
   \hline
  50 & BIAS & -0.20 & -0.37 & -0.38 & -0.31 & -0.42 & -0.07 & -0.18 & -0.16 & -0.12 & -0.16 \\
 & SD & 0.85 & 0.48 & 0.59 & 0.60 & 0.74 & 0.36 & 0.20 & 0.24 & 0.30 & 0.42 \\
 & RMSE & 0.88 & 0.60 & 0.70 & 0.67 & 0.85 & 0.37 & 0.27 & 0.29 & 0.32 & 0.45 \\
   \hline
  100 & BIAS & -0.07 & -0.16 & -0.13 & -0.31 & -0.28 & -0.04 & -0.08 & -0.08 & -0.18 & -0.07 \\
 & SD & 0.62 & 0.41 & 0.49 & 0.47 & 0.56 & 0.27 & 0.15 & 0.19 & 0.19 & 0.31 \\
 & RMSE & 0.62 & 0.44 & 0.51 & 0.56 & 0.62 & 0.27 & 0.17 & 0.21 & 0.26 & 0.32 \\
\hline
\multicolumn{12}{c}{OLS \& Random-$X$ Setting}\\
   \hline
25 & BIAS & -0.17 & -0.36 & -0.37 & -0.42 & -0.53 & -0.07 & -0.17 & -0.17 & -0.20 & -0.21 \\
 & SD & 1.16 & 0.67 & 0.83 & 0.94 & 0.97 & 0.53 & 0.28 & 0.35 & 0.44 & 0.56 \\
 & RMSE & 1.17 & 0.76 & 0.91 & 1.03 & 1.10 & 0.53 & 0.33 & 0.39 & 0.49 & 0.60 \\
   \hline
  50 & BIAS & -0.15 & -0.30 & -0.29 & -0.33 & -0.37 & -0.07 & -0.13 & -0.11 & -0.16 & -0.14 \\
 & SD & 0.87 & 0.50 & 0.65 & 0.67 & 0.73 & 0.39 & 0.20 & 0.27 & 0.33 & 0.41 \\
 & RMSE & 0.89 & 0.58 & 0.71 & 0.74 & 0.82 & 0.39 & 0.24 & 0.29 & 0.37 & 0.44 \\
   \hline
  100 & BIAS & -0.10 & -0.19 & -0.24 & -0.22 & -0.27 & -0.04 & -0.10 & -0.09 & -0.09 & -0.14 \\
 & SD & 0.59 & 0.37 & 0.46 & 0.54 & 0.55 & 0.26 & 0.15 & 0.19 & 0.24 & 0.30 \\
& RMSE & 0.60 & 0.42 & 0.52 & 0.58 & 0.61 & 0.27 & 0.18 & 0.21 & 0.25 & 0.33 \\
\hline

\end{tabular}}
\end{center}

\end{table}

\begin{table}[htbp!]

\caption{The average bias (BIAS),  standard deviation (SD), and  root mean
squared error (RMSE) of the constrained GLS and FGLS estimates obtained from Model A with  the mixture of normal random errors under the Fixed-$X$ and Random-$X$ settings.}
\label{tb:EC1m-2}
\begin{center}
\scalebox{0.99}{
\begin{tabular}{c|c|rrrrr|rrrrr}
  \hline
  & & \multicolumn{5}{c|}{$p=5$}     &    \multicolumn{5}{c}{$p=25$} \\
\cline{3-7}\cline{8-12}
 $n$  &\multicolumn{1}{c|}{}& \multicolumn{1}{c}{$\hat \beta_{1}$} & \multicolumn{1}{c}{$\hat \beta_{2}$} & \multicolumn{1}{c}{$\hat \beta_{3}$} & \multicolumn{1}{c}{$\hat \beta_{4}$} & \multicolumn{1}{c|}{$\hat \beta_{5}$} & \multicolumn{1}{c}{$\hat \beta_{1}$} & \multicolumn{1}{c}{$\hat \beta_{2}$} & \multicolumn{1}{c}{$\hat \beta_{3}$} & \multicolumn{1}{c}{$\hat \beta_{4}$} & \multicolumn{1}{c}{$\hat \beta_{5}$}   \\
   \hline

   \multicolumn{12}{c}{GLS \& Fixed-$X$ Setting}\\
   \hline
25 & BIAS & 0.06 & -0.03 & -0.00 & 0.00 & -0.02 & 0.02 & 0.02 & -0.01 & -0.03 & -0.04 \\
& SD & 1.30 & 0.62 & 0.64 & 0.64 & 0.73 & 0.55 & 0.21 & 0.27 & 0.29 & 0.44 \\
& RMSE & 1.30 & 0.62 & 0.64 & 0.64 & 0.73 & 0.55 & 0.21 & 0.27 & 0.29 & 0.44 \\
   \hline
  50 & BIAS & 0.02 & 0.02 & -0.00 & -0.03 & -0.03 & -0.03 & -0.02 & -0.01 & -0.01 & -0.02 \\
 & SD & 0.87 & 0.35 & 0.38 & 0.36 & 0.49 & 0.37 & 0.13 & 0.13 & 0.19 & 0.25 \\
 & RMSE & 0.87 & 0.35 & 0.38 & 0.36 & 0.49 & 0.37 & 0.13 & 0.13 & 0.19 & 0.25 \\
   \hline
  100 & BIAS & 0.02 & 0.02 & 0.02 & -0.05 & -0.02 & -0.00 & 0.00 & -0.00 & -0.00 & -0.00 \\
 & SD & 0.61 & 0.28 & 0.34 & 0.30 & 0.37 & 0.27 & 0.09 & 0.09 & 0.08 & 0.14 \\
 & RMSE & 0.61 & 0.29 & 0.34 & 0.30 & 0.37 & 0.27 & 0.09 & 0.09 & 0.08 & 0.14 \\
\hline
\multicolumn{12}{c}{GLS \& Random-$X$ Setting}\\
   \hline
25 & BIAS & -0.02 & -0.03 & -0.05 & -0.03 & -0.05 & -0.02 & -0.01 & -0.02 & -0.03 & -0.02 \\
 & SD & 1.28 & 0.56 & 0.64 & 0.68 & 0.75 & 0.53 & 0.21 & 0.25 & 0.29 & 0.39 \\
 & RMSE & 1.28 & 0.56 & 0.64 & 0.68 & 0.75 & 0.53 & 0.21 & 0.25 & 0.29 & 0.39 \\
   \hline
  50 & BIAS & -0.04 & -0.02 & -0.03 & -0.02 & -0.06 & -0.01 & -0.02 & 0.01 & -0.01 & -0.02 \\
 & SD & 0.86 & 0.37 & 0.46 & 0.46 & 0.50 & 0.38 & 0.14 & 0.17 & 0.21 & 0.26 \\
 & RMSE & 0.86 & 0.37 & 0.46 & 0.46 & 0.50 & 0.38 & 0.14 & 0.17 & 0.21 & 0.26 \\
   \hline
  100 & BIAS & -0.01 & -0.02 & -0.01 & -0.02 & -0.01 & -0.01 & -0.00 & -0.01 & -0.00 & -0.02 \\
 & SD & 0.60 & 0.25 & 0.31 & 0.33 & 0.33 & 0.27 & 0.10 & 0.12 & 0.15 & 0.18 \\
 & RMSE & 0.60 & 0.25 & 0.31 & 0.33 & 0.33 & 0.27 & 0.10 & 0.12 & 0.15 & 0.18 \\
\hline

   \multicolumn{12}{c}{FGLS \& Fixed-$X$ Setting}\\
   \hline
25 & BIAS & -0.11 & -0.29 & -0.27 & -0.34 & -0.35 & -0.07 & -0.08 & -0.09 & -0.12 & -0.17 \\
 & SD & 1.26 & 0.84 & 1.04 & 0.92 & 1.05 & 0.54 & 0.25 & 0.32 & 0.39 & 0.63 \\
 & RMSE & 1.26 & 0.89 & 1.08 & 0.98 & 1.10 & 0.55 & 0.26 & 0.33 & 0.41 & 0.65 \\
   \hline
  50 & BIAS & -0.13 & -0.16 & -0.18 & -0.17 & -0.32 & -0.03 & -0.04 & -0.06 & -0.05 & -0.09 \\
 & SD & 0.85 & 0.53 & 0.67 & 0.68 & 0.76 & 0.38 & 0.17 & 0.22 & 0.26 & 0.39 \\
 & RMSE & 0.86 & 0.55 & 0.69 & 0.70 & 0.83 & 0.38 & 0.18 & 0.23 & 0.26 & 0.40 \\
   \hline
  100 & BIAS & -0.09 & -0.09 & -0.11 & -0.14 & -0.13 & -0.04 & -0.03 & -0.04 & -0.05 & -0.06 \\
 & SD & 0.61 & 0.37 & 0.51 & 0.50 & 0.55 & 0.27 & 0.13 & 0.18 & 0.24 & 0.22 \\
 & RMSE & 0.61 & 0.38 & 0.52 & 0.52 & 0.56 & 0.28 & 0.13 & 0.18 & 0.24 & 0.23 \\
\hline
\multicolumn{12}{c}{FGLS \& Random-$X$ Setting}\\
   \hline
25 & BIAS & -0.22 & -0.29 & -0.32 & -0.34 & -0.37 & -0.06 & -0.08 & -0.10 & -0.10 & -0.12 \\
 & SD & 1.24 & 0.74 & 0.87 & 1.04 & 1.24 & 0.52 & 0.25 & 0.32 & 0.40 & 0.55 \\
 & RMSE & 1.26 & 0.79 & 0.93 & 1.09 & 1.30 & 0.53 & 0.27 & 0.33 & 0.41 & 0.57 \\
   \hline
  50 & BIAS & -0.11 & -0.14 & -0.18 & -0.21 & -0.23 & -0.05 & -0.05 & -0.04 & -0.06 & -0.09 \\
 & SD & 0.88 & 0.48 & 0.68 & 0.71 & 0.74 & 0.40 & 0.16 & 0.21 & 0.27 & 0.35 \\
 & RMSE & 0.89 & 0.50 & 0.70 & 0.74 & 0.78 & 0.40 & 0.17 & 0.22 & 0.28 & 0.36 \\
   \hline
  100 & BIAS & -0.07 & -0.10 & -0.11 & -0.13 & -0.17 & -0.04 & -0.03 & -0.04 & -0.03 & -0.05 \\
 & SD & 0.63 & 0.36 & 0.56 & 0.58 & 0.84 & 0.27 & 0.13 & 0.19 & 0.22 & 0.30 \\
 & RMSE & 0.63 & 0.37 & 0.57 & 0.60 & 0.85 & 0.27 & 0.13 & 0.19 & 0.22 & 0.30 \\
\hline

\end{tabular}}
\end{center}

\end{table}

Under Model A with random errors generated from the mixture of normal   distributions,  Tables \ref{tb:EC1m-1} -- \ref{tb:EC1m-2} present the average bias (BIAS), standard deviation (SD), and root mean squared error (RMSE) of the four  constrained  estimators via 1,000 realizations.
In addition, Table \ref{tb:EC1m-3} reports  two types averaged estimation errors in the covariance matrix estimate $C_{\hat\beta}(X_i)$. They are the S-Error ($\max_{1\leq i\leq n}\|C_{\hat\beta}(X_i)-C_{\beta_0}(X_i)\|_2$ with the spectral norm $\|\cdot\|_2$) and the F-error ($p^{-1/2}\max_{1\leq i\leq n}\|C_{\hat\beta}(X_i)-C_{\beta_0}(X_i)\|_F$  with the Frobenius norm $\|\cdot\|_F$), where $\hat\beta$ can be any one of the four estimates,
$\hat\beta_{\QMLE}^c$, $\hat\beta_{\OLS}^c$, $\hat\beta_{\GLS}^c$, and $\hat\beta_{\FGLS}^c$.
 Moreover, Tables \ref{tb:S1} -- \ref{tb:S3}  present the results  when random errors are generated from the standard normal distribution.

Tables \ref{tb:EC1m-1}, \ref{tb:EC1m-2} and \ref{tb:EC1m-3}
and Tables \ref{tb:S1} -- \ref{tb:S3} reveal five important findings. (I) The BIAS, SD, RMSE, S-Error and F-Error generally decrease as $n$ or $p$ gets larger, which confirms the consistency of the parameter estimators demonstrated in Theorems \ref{tm:2} -- \ref{tm:3}, and the consistency of  the covariance matrix estimators shown in Theorems \ref{tm:region} -- \ref{tm:cons}. It is worth noting that the biases of  the constrained QMLE, OLS and FGLS estimates  generally decrease as $n$ or $p$ gets larger, but they cannot be zero.   According to Remark \ref{re:positived}, this finding is reasonable since  the constrained estimators are asymptotically equivalent to the unconstrained estimators. Hence, the constrained estimators
are not unbiased. Rather, they are asymptotically unbiased. In this example, however,  the BIAS of the constrained GLS estimate  is very close to zero across all settings of  $n$ and $p$.
(II)  The SD and RMSE of the GLS and FGLS estimates get closer to  those of QMLE as either $n$ or $p$ increases, while these three estimates are generally superior to the OLS estimates. (III) The SD and RMSE of the GLS estimates are generally closer and comparable to  those of QMLE under
the normal random errors, while they are slightly inferior to
 those of QMLE under
the mixture of normal random errors.   (IV) The BIAS, SD, RMSE, S-Error and F-Error of the FGLS estimate gradually approach those of the GLS estimate as both $n$ and $p$ get larger, which is consistent with the comment below Theorem \ref{tm:3}.
 It is of interest to note that the GLS estimates sometimes have higher RMSE than that of the OLS or FGLS estimates, especially when $n$ and $p$ are small. One possible explanation is due to the Monte Carlo approximation errors. Another  possibility is that some errors may be due to small $n$ and $p$. This is because
the GLS is only asymptotically as efficient as the FGLS, and it is more efficient than the OLS as $n$ or $p$ goes to infinity. (V) There is not much difference
between the Fixed-$X$ and Random-$X$ settings in the RMSE, S-Error and F-Error measures across all four estimates. This finding is not surprising  since Theorems \ref{tm:2} -- \ref{tm:cons} hold for both fixed and random settings.

After analyzing the estimation performance for Model A, we  examine  the performance of the parameter estimates under Model B. Tables \ref{tb:S4} -- \ref{tb:S6} report the results  when the random errors are generated from the mixture of normal distributions, while Tables \ref{tb:S7} -- \ref{tb:S9} present the results  when the random errors are simulated from the standard normal distribution. These tables show qualitatively similar findings to those in  Tables \ref{tb:EC1m-1}, \ref{tb:EC1m-2} and \ref{tb:EC1m-3}
and Tables \ref{tb:S1} -- \ref{tb:S3}.

In addition to parameter estimation, we next compare 
the expected test error of covariance regression models under the Fixed-$X$ and Random-$X$ settings. The expected test error is  empirically calculated by averaging the test errors (\ref{eq:teste}) from 1,000 realizations, which is the empirical version of $\Err$ defined in (\ref{eq:prerr}).
Under   the   mixture of normal  random errors,
 Figure  \ref{fig:PC1m-1} depicts the corresponding expected test errors.
Note that
 Figure \ref{fig:PC1m-1} is obtained under Model A, and  
the  expected test error is calculated by setting   $\bmW(X_\alpha)= \bmV^-(X_\alpha)$ for $\alpha=1,\cdots,n$ and $01,\cdots,0n$, and $\hat \bmC(\cdot)=\bmC_{\hat\beta_{\GLS}^c}(\cdot)$ in (\ref{eq:teste}).
To assess the impact of model misspecification under Model A,  we consider  the correctly-specified model $C_{\beta}(X_\alpha)=\beta^{(1)}I_p+\sum_{k=2}^{5}\beta^{(k)}X_\alpha^{(k)}$ and a misspecified model $C_{\beta}(X_\alpha)=\beta^{(1)}I_p+\sum_{k=2}^{4}\beta^{(k)}X_\alpha^{(k)}+\beta^{(5)}\overline X_\alpha^{(5)}$,  where $\overline X_\alpha^{(5)}=(\overline X_{\alpha,j_1j_2}^{(5)})_{p\times p}$ with $\overline X_{\alpha,j_1j_2}^{(5)}=(X_{\alpha,j_1j_2}^{(5)})^2$. In addition,  Figure \ref{fig:PC1n-1} presents the results  under the standard normal random errors.

For $p=5$, Figure \ref{fig:PC1m-1} and Figure \ref{fig:PC1n-1} reveal the following important findings.
(I) The Random-$X$ expected test error,  $\ErrR$, is larger than the Fixed-$X$ expected test
error, $\ErrF$, in both the  mixture of  normal and the normal random errors.
 Note that $\ErrF$ under our Fixed-$X$ simulation setting converges to $\ErrS$ based on   Proposition \ref{pn:Err} in Section \ref{sec:adS3} of the supplementary material. Hence, this finding is consistent with the theoretical result in Theorem \ref{eq:excess}. In addition,
the magnitude of the difference between  $\ErrR$ and $\ErrF$  under the  mixture of normal random errors is greater than that under normal random errors.
(II) The difference between  $\ErrR$ and $\ErrF$ becomes small as $n$ increases when the model is correctly specified. This finding supports  Remark \ref{re:5}. However, the difference
increases  when the model is misspecified;  see  the ``Remark Regarding Theorem \ref{eq:excess}'' in Section \ref{sec:adS3} of the supplementary material.
It is worth noting that
Figure \ref{fig:PC1m-1} and Figure \ref{fig:PC1n-1} indicate that
the aforementioned difference between  $\ErrR$ and $\ErrF$ is not as strong
for  $p=25$ as for $p=5$.
 However, there still exist  differences  between $\ErrR$ and $\ErrF$ when the model is misspecified. These findings are consistent with  the ``Remark Regarding Theorem \ref{eq:excess}'' in Section \ref{sec:adS3}.

 Subsequently, we evaluate the performance of  the estimates of expected  test errors,
$\Cp$, $\widehat\Cp$, $\RCp$, $\widehat\RCp$ and $\OCV$, when the data are generated under the Random-$X$ setting.
As mentioned in Section \ref{sec:non-negative},  $\Cp$ and $\widehat\Cp$ are
obtained under Fixed-$X$, while  $\RCp$ and $\widehat\RCp$ are obtained under Random-$X$.
Accordingly, mis-using  $\Cp$ and $\widehat\Cp$ under Random-$X$ can underestimate the expected test error compared to that of $\RCp$ and $\widehat\RCp$.
In addition, we compute the average values of  $\Cp$, $\widehat\Cp$, $\RCp$, $\widehat\RCp$, and  $\OCV$ from 1,000 realizations, and then subtract  the
corresponding expected test errors.  Accordingly, we obtain
empirical  biases that are a good measure for assessing  performance.
Figure \ref{fig:PC1m-2} reports the  biases of $\Cp$, $\widehat\Cp$, $\RCp$, $\widehat\RCp$ and $\OCV$ in absolute value terms for  the  mixture of normal random errors. In addition,  Figure \ref{fig:PC1n-2} presents the results  for the standard normal random errors.
These two figures reveal three important findings.  (I) The absolute bias of $\RCp$ ($\widehat\RCp$) is smaller than that of $\Cp$ ($\widehat\Cp$) across all settings. This finding is not surprising since $\Cp$ and $\widehat\Cp$ underestimate expected test errors under the Random-$X$  setting; see the comment below equation (\ref{pn:CVRE}).
 (II) $\RCp$ and $\widehat\RCp$ generally have a smaller absolute bias  compared
to $\OCV$, though as $n$ gets larger, the bias of $\OCV$ decreases.   This finding is reasonable since $\OCV$ is designed  for the general covariance regression model, while $\RCp$ and $\widehat\RCp$ are obtained for the linear covariance regression model, which is considered in this simulation study. (III) The absolute bias of $\widehat\RCp$ generally gets closer to $\RCp$ as $n$ gets larger, which is sensible since the model is correctly specified; see the comment before equation (\ref{pn:CVRE}) and in Remark \ref{re:4}.  Based on these findings, we recommend using $\widehat\RCp$ in practice.

The above studies are based on Model A with the setting $\bmW(X_\alpha)= \bmV^-(X_\alpha)$ and the estimate $\hat \bmC(\cdot)=\bmC_{\hat\beta_{\GLS}^c}(\cdot)$. We next consider the setting $\bmW(X_\alpha)= I_p$ with the estimate  $\hat \bmC(\cdot)=\bmC_{\hat\beta_{\OLS}^c}(\cdot)$ under the same model.
The comparisons between the expected test errors  $\ErrR$ and $\ErrF$ as
well as the examinations of the absolute biases of  $\Cp$, $\widehat\Cp$, $\RCp$, $\widehat\RCp$ and $\OCV$, are presented in Figures \ref{fig:S1} -- \ref{fig:S3} and
Figures \ref{fig:S51} -- \ref{fig:S4} for  the mixture of normal random errors and  the standard normal random errors, respectively. All of these figures   are qualitatively similar to Figures
\ref{fig:PC1m-1} -- \ref{fig:PC1m-2} and Figures \ref{fig:PC1n-1} -- \ref{fig:PC1n-2}, although some of the differences between $\ErrR$ and $\ErrF$ and  among $\Cp$, $\widehat\Cp$, $\RCp$, $\widehat\RCp$ and $\OCV$ are not as strong as   depicted in Figures
\ref{fig:PC1m-1} -- \ref{fig:PC1m-2} and Figures \ref{fig:PC1n-1} -- \ref{fig:PC1n-2}.

Finally, we compare  expected test errors based on Model B under the following two different settings:
(i)  $\bmW(X_\alpha)= \bmV^-(X_\alpha)$ with the estimate $\hat \bmC(\cdot)=\bmC_{\hat\beta_{\GLS}^c}(\cdot)$; and
(ii) $\bmW(X_\alpha)= I_p$ with the estimate $\hat \bmC(\cdot)=\bmC_{\hat\beta_{\OLS}^c}(\cdot)$.
 The  comparisons between the  expected test errors  $\ErrR$ and $\ErrF$ as well as the examinations of the absolute biases
are presented in Figures \ref{fig:S5} -- \ref{fig:S11b} under setting (i), and in Figures \ref{fig:S11} -- \ref{fig:S12} under setting (ii). Compared to Figures
\ref{fig:PC1m-1} -- \ref{fig:PC1n-2}, the differences between $\ErrR$ and $\ErrF$ and  among $\Cp$, $\widehat\Cp$, $\RCp$, $\widehat\RCp$ and $\OCV$ in Figures \ref{fig:S5} -- \ref{fig:S11b} are less pronounced. On the other hand, the performances of Figures  \ref{fig:S11} -- \ref{fig:S12} are qualitatively similar to those of Figures \ref{fig:S1}  -- \ref{fig:S4}.  In addition to these results, additional simulation studies according to anonymous referees' suggestions are presented in Section \ref{sec:Ssimu} of the supplementary material.

\section{Conclusion}

\label{sec:5}

 In this article, we study two challenging tasks in covariance regression models. One is establishing asymptotic theories for parameter estimation in random covariance regression models; and the other is building up model assessment theories for assessing differences between the fixed and random covariance regression model settings.
Specifically, under the Random-$X$ setting,
we demonstrate that both QMLE and WLS estimators are consistent and asymptotically normal as long as either $n$ or $p$ goes to infinity. As suggested by anonymous reviewers, we further obtain the asymptotic properties of the QMLE and WLS 
estimator  by  allowing several possibilities: a general error form for the data generating process, the
dimension  of $\beta$ (i.e., $K$) to go to infinity,  the model to be misspecified, and the conditional dependence of $Y_i$.  
In addition, the model assessment theory indicates that the difference between the expected test errors under the Fixed-$X$ and Random-$X$ settings is largely due to the bias and variance induced by the randomness in $X$. This finding allows us to propose two estimators, $\widehat \RCp$ and $\OCV$, to estimate 
Random-$X$ expected test errors.
 Simulation studies in Section \ref{sec:40} support both estimation  and 
model assessment theories.

The pioneering work of this paper opens a new avenue in
studying covariance regression models.
For example, this framework for estimation and model assessment in the random covariance regression model  (\ref{eq:model}) can be extended to analyze multivariate volatility models   $\Cov[Y_i|\mG_{i-1}]= C_\beta(X_{i-1})$ (see, e.g.,  Chapter 7 of \citealp{Tsay:2014}), where  $\mG_{i-1}$ is the  $\sigma$-algebra generated from the observations up to
time  $i-1$, and $X_i$ is a stochastic process.
Specifically, we can estimate $\Cov[Y_i|\mG_{i-1}]$ and then construct time-varying portfolios for asset returns, which is  an important subject in financial risk management; see, e.g., \citet{Tse2002}, \citet{Laurent2012}, 
and  \citet{Poignard2022}. 
 In addition, analogous estimation and model assessment approaches can be applied to mean-covariance regression models (see, e.g., \citealp{Anderson1973} and \citealp{Zou:Lan:Li:Tsai2021}). 
For instance, we can adapt the  mean-covariance regression model from  \citet{Zou:Lan:Li:Tsai2021}, and consider $Y_{\alpha}=BZ_{\alpha}+V_{\alpha}$, where $Z_{\alpha}\in\mathbb{R}^d$ is the $d$-dimensional explanatory vector, $B$ is the $p\times d$ mean regression coefficient matrix, and $V_{\alpha}$ is standardized to have conditional mean zero and covariance matrix $\Cov[V_{\alpha}|\mF_{\mX}]= C_\beta(X_{\alpha})$. Furthermore, 
$\big\{(X_\alpha,Z_\alpha,Y_{\alpha}):\alpha \in\{1,\cdots,n\}\big\}$ and 
$\big\{(X_\alpha,Z_{\alpha},Y_{\alpha}):\alpha \in\{01,\cdots,0m\}\big\}$ are 
training data and test data, respectively, and 
$\mF_{\mX}=\sigma\langle(X_1,Z_1),\cdots,(X_n,Z_n),(X_{01},Z_{01}),\cdots, (X_{0m},Z_{0m})\rangle$. Accordingly, $\rmE[Y_{0i}|\mF_\mX]=BZ_{0i}$.
Then, employing  equation (3.2) of \citet[p. 321]{Zou:Lan:Li:Tsai2021}, we  obtain the MLE of $\rmvec (B)$ below,
\[
\rmvec (\hat B)= \lbk  \sum_{i=1}^n (Z_i\otimes I_p) C_{\hat\beta}^{-1}(X_i) (Z_i^\top\otimes I_p)\rbk^{-1}\sum_{i=1}^n (Z_i\otimes I_p) C_{\hat\beta}^{-1}(X_i) Y_i,
\]
where $\hat\beta$ is the estimator of covariance regression parameter $\beta$.
As a result, the prediction of $Y_{0i}$ is $\hat BZ_{0i}$, for $i=1,\cdots,m$.
It is clear that  the prediction $\hat BZ_{0i}$ depends on the covariance regression structure $C_{\hat\beta}(\cdot) $. Hence, utilizing the estimation and model assessment approaches developed  in this paper for covariance regression,
one can   select better models for $C_\beta(\cdot)$ and make more accurate 
inferences. In sum, a reliable $C_{\hat\beta}(\cdot) $ can play an important role in improving the accuracy of the  mean  prediction $\hat BZ_{0i}$. Analogous mean predictions can be considered for linear mixed models (see, e.g., \citealp{McCulloch2008} and \citealp{Demidenko2013}).

It is also worth noting that  in the context of nonparametric statistics,  machine learning, and econometrics,  theories are commonly established under the Random-$X$ setting (see, e.g., \citealp{Gyorfy2002}, \citealp{Vapnik1998} and \citealp{Wooldridge2002}).  Accordingly, the proposed theoretical framework can facilitate extending covariance regression to nonparametric, machine learning, and econometric models for future research. Specifically, we can generalize the nonparametric variance function estimation (see, e.g., \citealp{Wang2008}) to nonparametric covariance function estimation, via nonparametric covariance regression models.  Finally,  \citeauthor{Shen2020}'s (\citeyear{Shen2020}) two-level parametric bootstrap approach can be employed to estimate the bias difference $\mB^{(d)}$ and variance difference $\mV^{(d)}$ for covariance regression models. Compared to $\widehat \RCp$, this type of approach does not require the model to be correctly specified, and the estimator of $\beta$ can be different from the WLS estimator.
We believe  these extensions would strengthen the usefulness of the random covariance regression model.

\section*{References}

\begin{description}
\newcommand{\enquote}[1]{``#1''}
\expandafter\ifx\csname natexlab\endcsname\relax\def\natexlab#1{#1}\fi

\bibitem[{Anderson, 1973}]{Anderson1973}
Anderson, T. W. (1973).
\enquote{Asymptotically efficient estimation of covariance matrices with linear structure,}
\textit{The Annals of Statistics}, 1, 135-141.

\bibitem[{Bai and Silverstein, 2004}]{Bai2004}
Bai, Z. D. and Silverstein, J. W. (2004).
\enquote{Central limit theorem for linear spectral statistics of large dimensional sample covariance matrices,}
\textit{The Annals of Probability}, 32, 553-605.

\bibitem[{Bollerslev and Wooldridge, 1992}]{Bollerslev1992}
Bollerslev, T. and Wooldridge, J. M. (1992).
\enquote{Quasi-maximum likelihood estimation and inference in dynamic models with time-varying covariance,}
\textit {Econometric Reviews}, 11, 143-173.

\bibitem[{Breiman and Spector, 1992}]{Breiman1992}
Breiman, L. and Spector, P. (1992).
\enquote{Submodel selection and evaluation in regression. The X-random case,}
\textit{International statistical review/revue internationale de Statistique}, 60, 291-319.

\bibitem[{Cai et al.(2020)}]{Cai2020}
Cai, T. T., Hu, J., Li, Y. and Zheng, X. (2020).
\enquote{High-dimensional minimum variance portfolio estimation based on high-frequency data,}
\textit{Journal of Econometrics,} 214, 482-494.

\bibitem[{Chiu et al.(1996)}]{Chiu1996}
Chiu, T. Y., Leonard, T. and Tsui, K. W. (1996).
\enquote{The matrix-logarithmic covariance model,}
\textit{Journal of the American Statistical Association,} 91, 198-210.

\bibitem[{Chung(2001)}]{Chung2001}
Chung, K. L. (2001).
\textit{A Course in Probability Theory}, Academic Press.

\bibitem[{Cressie(1993)}]{Cressie1993}
Cressie, N. A. C. (1993).
\textit{Statistics for Spatial Data}, New York: Wiley.

\bibitem[{Demidenko(2013)}]{Demidenko2013}
Demidenko, E. (2013).
\textit{Mixed Models: Theory and Applications with R}, John Wiley \& Sons.

\bibitem[{DeMiguel et al.(2019)}]{DeMiguel2009}
DeMiguel, V., Garlappi, L. and Uppal, R. (2009). \enquote{Optimal versus naive diversification: How inefficient is the 1/N portfolio strategy?,}
\textit{The Review of Financial Studies},  22, 1915-1953.

\bibitem[{Fan et al.(2012)}]{Fan2012}
Fan, J., Li, Y. and Yu, K. (2012). \enquote{Vast volatility matrix estimation using high-frequency data for portfolio selection,}
\textit{Journal of the American Statistical Association},  107, 412-428.

\bibitem[{Fox and Dunson(2015)}]{Fox2015}
Fox, E. B. and Dunson, D. B. (2015).  \enquote{Bayesian nonparametric covariance regression,}
\textit{The Journal of Machine Learning Research}, 16, 2501-2542.

\bibitem[{Gy\"orfy et al.(2002)}]{Gyorfy2002}
Gy\"orfy, L., Kohler, M., Krzyzak, A. and Walk, H.  (2002).
\textit{A Distribution-Free Theory of Nonparametric Regression}, New York: Springer.

\bibitem[{Hastie et al.(2009)}]{Hastie2009}
Hastie, T., Tibshirani, R. and Friedman, J. (2009).
\textit{The Elements of Statistical Learning: Data Mining, Inference and Prediction}, 2nd Edition, New York: Springer.

\bibitem[{Hoff and Niu(2012)}]{Hoff2012}
Hoff, P. D. and Niu, X. (2012). \enquote{A covariance regression model,}
\textit{Statistica Sinica}, 22, 729-753.

\bibitem[{Huang et al.(2019)}]{Huang2019}
Huang, D., Lan, W., Zhang, H. H. and Wang, H. (2019). \enquote{Least squares estimation of spatial autoregressive models for large-scale social networks,}
\textit{Electronic Journal of Statistics}, 13, 1135-1165.

\bibitem[{Kelejian and Prucha(2001)}]{Kelejian2001}
Kelejian, H. and Prucha, I. (2001) \enquote{On the asymptotic distribution of the Moran I test statistic with applications,} \textit{Journal of Econometrics}, 104, 219-257.

\bibitem[{Laurent et al.(2012)}]{Laurent2012}
Laurent, S., Rombouts, J. V. and Violante, F. (2012). \enquote{On the forecasting accuracy of multivariate GARCH models,}
\textit{Journal of Applied Econometrics}, 27, 934-955.

\bibitem[{Lee, 2004}]{Lee2004}
Lee, L. F. (2004).
\enquote{Asymptotic distributions of quasi-maximum likelihood estimators for spatial autoregressive models,}
\textit{Econometrica}, 72, 1899-1925.

\bibitem[{LeSage and Pace(2009)}]{LeSage2009}
LeSage, J. and Pace, R. K.  (2009).
\textit{Introduction to Spatial Econometrics}, New York: Chapman \& Hall.

\bibitem[{Ma et al.(2019)}]{Ma2019}
Ma, Y., Pan, R., Zou, T. and Wang, H. (2019). \enquote{A naive least squares method for spatial autoregression with covariates,}
\textit{Statistica Sinica},  30, 653-672.

\bibitem[{Mallows, 1973}]{Mallows1973}
Mallows, C. L. (1973).
\enquote{Some comments on $C_p$,}
\textit{Technometrics}, 15, 661-675.

\bibitem[{Markowitz(1952)}]{Markowitz1952}
Markowitz, H. (1952). \enquote{Portfolio selection,}
\textit{The Journal of Finance}, 7, 77-91.

\bibitem[{McCulloch et al.(2008)}]{McCulloch2008}
McCulloch, C. E., Searle, S.R. and Neuhaus, J. M.  (2008).
\textit{Generalized, Linear, and Mixed Models}, 2nd Edition, John Wiley \& Sons.

\bibitem[{Pan and Zhou, 2008}]{Pan2008}
Pan, G. and Zhou, W. (2008).
\enquote{Central limit theorem for signal-to-interference ratio of
reduced rank linear receiver,}
\textit{The Annals of Applied Probability}, 18, 1232-1270.

\bibitem[{Poignard and Asai, 2023}]{Poignard2022}
Poignard, B. and Asai, M. (2023).
\enquote{High‐dimensional sparse multivariate stochastic volatility models,}
\textit{Journal of Time Series Analysis}, 44, 4-22.

\bibitem[{Rosset and Tibshirani, 2020}]{Rosset2020}
Rosset, S. and Tibshirani, R. J. (2020).
\enquote{From Fixed-$X$ to Random-$X$ regression: Bias-variance decompositions, covariance penalties, and prediction error estimation,}
\textit{Journal of the American Statistical Association}, 115, 138-151.

\bibitem[{Scott and Handcock, 2001}]{Scott2001}
Scott, M. A. and Handcock, M. S. (2001).
\enquote{Covariance models for latent structure in longitudinal data,}
\textit{Sociological Methodology}, 31, 265-303.

\bibitem[{Shen and Huang, 2020}]{Shen2020}
Shen, X. and Huang, H. C. (2020).
\enquote{Discussion of `From Fixed-X to Random-X Regression: Bias-Variance Decompositions, Covariance Penalties, and Prediction Error Estimation',}
\textit{Journal of the American Statistical Association}, 115, 152-156.

\bibitem[{Tsay, 2014}]{Tsay:2014}
Tsay, R. S. (2014).
\newblock{\em Multivariate Time Series Analysis},
\newblock{Hoboken: John Wiley \& Sons.}

\bibitem[{Tse and Tsui, 2002}]{Tse2002}
Tse, Y. K. and Tsui, A. K. C. (2002).
\enquote{A multivariate generalized autoregressive conditional heteroscedasticity model with time-varying correlations,}
\textit{Journal of Business \& Economic Statistics}, 20, 351-362.

\bibitem[{van der Vaart(1998)}]{van1998}
van der Vaart, A. W. (1998).
\newblock{\em  Asymptotic Statistics},
\newblock{Cambridge University Press.}

\bibitem[{Vapnik(1998)}]{Vapnik1998}
Vapnik, N. V. (1998).
\textit{Statistical Learning Theory}, New York: Wiley.

\bibitem[{Wang et al.(2008)}]{Wang2008}
Wang, L., Brown, L. D., Cai, T. T. and Levine, M. (2008).
\enquote{Effect of mean on variance function estimation in nonparametric regression,}
\textit{The Annals of Statistics},  36, 646-664.

\bibitem[{Wikle et al.(2019)}]{Wikle2019}
Wikle, C. K., Zammit-Mangion, A. and Cressie, N. (2019).
\textit{Spatio-Temporal Statistics with R}, Chapman and Hall.

\bibitem[{White, 1982}]{White1982}
White, H. (1982).
\enquote{Maximum likelihood estimation of misspecified models,}
\textit{Econometrica}, 50, 1-25.

\bibitem[{Wooldridge(2002)}]{Wooldridge2002}
Wooldridge, J. M. (2002).
\textit{Econometric Analysis of Cross Section and Panel Data}, MIT Press.

\bibitem[{Wooldridge(2012)}]{Wooldridge2012}
Wooldridge, J. M. (2012).
\textit{Introductory Econometrics: A Modern Approach}, Mason, OH: Cengage Learning.

\bibitem[{Zheng, 2012}]{Zheng2012}
Zheng, S. (2012).
\enquote{Central limit theorems for linear spectral statistics of large dimensional F-matrices,}
\textit{Annales de l'IHP Probabilit\'es et statistiques}, 48, 444-476.

\bibitem[{Zhou et al.(2017)}]{Zhou2017}
Zhou, J., Tu, Y., Chen, Y. and Wang, H. (2017). \enquote{Estimating spatial autocorrelation with sampled network data,}
\textit{Journal of Business \& Economic Statistics}, 35, 130-138.

\bibitem[{Zou et al.(2022)}]{Zou:Lan:Li:Tsai2021}
Zou, T., Lan, W., Li, R. and Tsai, C.-L. (2022). \enquote{Inference on covariance-mean regression,}
\textit{Journal of Econometrics}, 230, 318-338.

\bibitem[{Zou et al.(2017)}]{Zou:Lan:Wang:Tsai:2017}
Zou, T., Lan, W., Wang, H. and Tsai, C.-L. (2017). \enquote{Covariance regression analysis,}
\textit{Journal of the American Statistical Association}, 112, 266-281.

\end{description}

\end{document}